\documentclass[a4paper]{article}
\usepackage[utf8x]{inputenc}
\usepackage{fullpage}
\usepackage{amsthm}
\usepackage{amssymb}
\usepackage{amsmath}
\usepackage{amsfonts}
\usepackage{latexsym}
\usepackage{url}
\usepackage{dutchcal}
\usepackage{tkz-tab}
\usepackage[linesnumbered,ruled,vlined]{algorithm2e}
\usepackage{algorithmic}
\usepackage{paralist} 
\usepackage{tikz}
\usetikzlibrary{calc,fit,spy}
\usetikzlibrary{matrix}
\usepackage{pgfplots}
\usetikzlibrary{arrows.meta}
\usepgfplotslibrary{fillbetween}
\pgfplotsset{compat=1.8}
\usepackage{multirow}
\usepackage{array}

 \newcommand\ForAuthors[1]
 {\par\smallskip                     
  \begin{center}
   \fbox
   {\parbox{0.9\linewidth}
    {\raggedright\sc--- #1}
   }
  \end{center}
  \par\smallskip                     %
 }
  \newcommand\comment[1]{}

\title{Quadratic Maximization of Reachable Values of Stable Discrete-Time Affine Systems}
\author{Assalé Adjé\\ LAboratoire de Modélisation et Pluridisciplinaire et Simulations\\
LAMPS\\ Université de Perpignan Via Domitia\\ France
}
\def\mm{\mathbcal{e}}
\def\inte{\operatorname{int}}
\def\conv{\operatorname{conv}}
\def\rr{\mathbb R}
\def\rd{\rr^d}
\def\nn{\mathbb N}
\def\xin{X^\mathrm{in}}
\def\rea {\mathcal{R}}
\def\klyap{\mathbf{K}^{\rm lyap}}
\def\kdiag{\mathbf{K}^{\rm diag}}
\def\vdiag{V_{\rm diag}}
\def\lyap#1{\mathcal{L}_{#1}}

\def\ptset{\mathcal{L}(A,Q)}

\def\norm#1{\| #1\|}
\def\qex{Q_{\rm ex}}
\def\ptex{\mathcal{L}(A_\gamma,\qex)}

\def\lyapex{\rm ExL}
\def\cont{\mathcal{T}_{Q}}
\def\lyopt{V^{\rm opt}}
\def\cc{\mathbb C}
\def\cd{\cc^d}
\def\cdd{\cc^{d\times d}}

\def\hdd{\mathcal{H}^{d\times d}}
\def\hp{\hdd_+}
\def\hddp{\hdd_{++}}

\def\setsdp{C_{\alpha,\beta}}

\def\lmax#1{\lambda_{\operatorname{max}}\left(#1\right)}
\def\lmin#1{\lambda_{\operatorname{min}}\left(#1\right)}
\def\mus#1{\mu\left(#1\right)}
\def\topt#1{t_Q^*\left(#1\right)}
\def\sdp#1{C_{#1}}
\def\betaf{\mathbcal{b}}
\def\alphaf{\mathbcal{a}}
\def\agmx{\operatorname{Argmax}}
\def\co{c_0(\rr)}
\def\gse{\Delta}
\def\gs{\Delta^s}
\def\bigkse{\mathbf{K}}
\def\bigks{\mathbf{K}^s}
\def\poseq{\mathrm{Pos}}
\def\pos{\mathrm{Pos}^{s}}
\def\ks{k^{s}}
\def\kse{k}
\newcommand{\Min}{\operatorname*{Min}}
\newcommand{\st}{\operatorname*{s. t.}}
\def\nuopt{\nu_{\rm opt}}
\def\kopt{k_{\rm opt}}
\def\xopt{x_{\rm opt}}

\def\bigklyap{\mathbf{K}^{\rm lyap}}
\def\Idd{\operatorname{Id}}
\def\btilde{\tilde{b}}

\def\pab{P_{a,b}}
\def\pubp{P_{1,b'}}
\def\pub{P_{1,b}}
\def\bP{\overline{P}}
\def\oB{\overline{B}}

\newcommand{ \pmax}{P_{\mathrm{max}}}
\newcommand{ \pmaxi}{P_{\mathrm{max}}^i}
\def\eigen{\mathrm{Eigen}}
\newenvironment{psmallmatrix}
  {\left(\begin{smallmatrix}}
  {\end{smallmatrix}\right)}
\newtheorem{assumption}{Assumption}
\newtheorem{proposition}{Proposition}
\newtheorem{theorem}{Theorem}
\newtheorem{corollary}{Corollary}
\newtheorem{problem}{Problem}
\newtheorem{lemma}{Lemma}
\newtheorem{remark}{Remark}
\newtheorem{example}{Example}
\SetKwInOut{Input}{Input}
\SetKwInOut{Output}{Output}
\SetKwInOut{init}{Initialisation}
\begin{document}
\maketitle

\begin{abstract}
In this paper, we solve a maximization problem where the objective function is quadratic and the constraints set is the reachable values set of a stable discrete-time affine system. This problem is equivalent to solve an infinite number of linearly constrained quadratic maximization programs. To solve exactly and in finite time this problem, we have to safely extract a finite number of them. Safely means that we must guarantee that the optimal solution can be found within this extracted family. This family has to be the smallest possible. Therefore, we construct an integer-valued function defined on the solutions of the discrete Lyapunov equation. Those integers represent overapproximations of the number of quadratic programs to solve to obtain an optimal solution of our specific maximization problem. The integer-valued function is minimized in order to get the smallest possible overapproximation. The method proposed in the paper is first illustrated on a class of nondiagonalizable systems of dimension two and finally experimented on randomly generated instances of the problem.
\end{abstract}

\section{Introduction}
In this paper, we are interested in maximizing the quadratic images of the reachable values of a time-invariant discrete-time affine system. Optimization over the reachable values of dynamical systems arises to analyze the robustness or the performance of the system~\cite{petersen,iwasaki};  for inverse problems in control~\cite{ZHANG2019108593}; to verify some safety properties for the system~\cite{fazlyab2020safety}. The optimal value in this context can be viewed as the value which penalizes the most the system with respect to some criteria. The specialization of our maximization problem to quadratic objective functions and affine discrete-time systems comes from diverse problems in control theory. More precisely, our maximization problem appears in the analysis of peaks~\cite{peakArt1,balakrishnan1992computing} in which the quadratic form is the square of Euclidean norm. In this paper, we develop a generic approach including any square of quadratic norms. Furthermore, in~\cite{peakArt1}, the authors only propose an upper-bound for their maximization problem. Our problem also arises to prove ellipsoids safety for the system (e.g.~\cite{ghosh2019robust,kong2018reachable,belta2017formal}) meaning that the reachable values set has to be proved to be included in some given ellipsoid. In this case, the objective function of the problem is \emph{convex}. Dually, our formulation can be used to prove ellipsoids avoidance meaning that no reachable values belong to a given ellipsoid. In this case, the objective function of the problem is \emph{concave}. The method developed in this paper also includes the case where the objective function is \emph{linear}. 

In this paper, some assumptions are made on the system. The system is supposed to be time-invariant and autonomous (uncontrolled) without inputs nor disturbances. Moreover, we assume that the initial conditions are only assumed to be in some given polytope. Hence, a reachable value is a term of a recurrent sequence starting from a vector in this initial polytopic set. The set of constraints of our maximization problem is, thus, not closed nor convex and cannot be represented in machine. The quadratic objective function is, on the other hand, completely general i.e. can be non-homogeneous, convex, concave, indefinite or even linear. 
This context does not allow to use directly classical quadratic maximization methods. Additional restrictive hypotheses are made in the paper. The most restrictive hypothesis is the stability of the system i.e. the matrix defining the dynamics has a spectral radius strictly smaller than one. Even with the stability condition, an optimal solution may not exist. The existence of an optimal solution can be ensured by the existence of a positive term of the sequence of the interest. However, the optimal value is always finite. The main theoretical difficulty of the problem is the fact that we have to compare an infinite number of values to obtain the optimal value of the maximization problem. We then develop a formula depending on the data of the problem in order to only compare a finite number of values and to solve exactly the problem. 

The author of the paper initiated a work on the computation of an upper bound of the optimal value of the problem using semidefinite programming without any guarantees on the exactness of the overapproximation. The approach has been developed when the discrete-time system was piecewise affine~\cite{10.1007/978-3-319-54292-8_2} or polynomial~\cite{adje2015property}. The technique developed in this paper is proved to provide the exact optimal value.
This current paper generalizes the previous one where the matrix defining the dynamics was supposed to be diagonalizable~\cite{DBLP:journals/jota/Adje21}. In this paper, we still assume the stability of the system. This stability assumption permits to exploit the existence of a solution for the discrete Lyapunov equation. No additional hypotheses are made on the matrix defining the system.  

The closest work seems to be the one proposed by Ahmadi and G\"unl\"uk~\cite{7403149,ahmadi2018robust}. The authors are interested in solving an optimization problem similar to the one presented here. They consider linear objective functions and linear systems (or belongs to a finite family of some linear systems).
The formulation differs from ours. First, in their context, the state variable of the system has to remain in a polyhedral invariant whereas, in our context, a constraint is only imposed to initial conditions. Second, their problem deals with linear objective function. Moreover, the authors propose the computation of inner and outer approximations of the reachable values set based on semidefinite programming. This is not mandatory here. Finally, their approach can be used for switched linear systems~\cite{sun2006switched}. The main similarity is the computation of an upper bound on the number of iterations required to solve the problem. Those bounds are not comparable with the one proposed here since our frameworks are different. Some hypotheses made in the paper (the existence of a positive term) are connected to some decision problems for discrete-time dynamical systems~\cite{fijalkow2019decidability}. Those decision problems (Skolem problem and its variants) are still open~\cite{ouaknine2014positivity}. We do not provide any result about the decidability of the existence of positive terms. The goal of this paper is to develop a constructive method to solve computationally an optimization problem for which a positive term of a sequence must exist.  

The main contribution of the paper is the extension of the previous paper~\cite{DBLP:journals/jota/Adje21} to general stable affine systems (to allow
the nondiagonalizable case). We prove that the formula obtained in this paper is equal, once restricted to diagonalizable case, to the formula of the previous paper~\cite{DBLP:journals/jota/Adje21}. Some refinements have been made about the analysis of the formula. We also introduce the problem of the minimization of the formula with respect to solutions of the discrete Lyapunov   
equation. We propose a theoretical answer to this new nonlinear semidefinite program. The results of the paper are also illustrated on a class of linear systems with nondiagonalizable matrix dynamics. Benchmark tables are also provided to show the correctness of the approach.

The paper is organized as follows. Section~\ref{background} introduces the problem, the main results obtained in the previous paper~\cite{DBLP:journals/jota/Adje21}, some useful results about the maximizer of sequences which converge to 0 and the running example used in the paper. Section~\ref{mainresults} describes the main results of the paper. More precisely, we present the extension of the formula and the analysis of the new formula. We also illustrate the results and the tools on the running example. Section~\ref{minlyapunov} proposes a study of the minimization problem of the formula with respect to the Lyapunov candidates. We prove that this minimization problem is equivalent to a certain nonlinear semidefinite program for which an optimal solution exists. Section~\ref{comparison} compares the restriction of the formula obtained in this paper to diagonalizable matrices and the one proposed in~\cite{DBLP:journals/jota/Adje21}. Section~\ref{benchexamples} illustrates and analyses the formula on some concrete example and on some randomly generated problems. Section~\ref{conclusion} concludes the paper and gives some future directions.

\section{Background}
\label{background}
\subsection{Problem Formulation}
The motivation of the paper is the resolution of a maximization problem for which the objective function is a quadratic form and the constraint set is the reachable values set of a stable (time-invariant) affine discrete-time system. The system is charatacterized by, a $d\times d$ matrix $A$, a $\rd$-vector $b$ and a polytope $\xin$ which is the set of initial conditions. We consider the following recurrence relation starting from $x_0\in\xin$ :
\begin{equation}
\label{stateeq}
  x_{k+1}=Ax_k+b,\ \forall\, k\in\nn \enspace.
\end{equation}
Then, we introduce the set $\rea_k$ of vectors that are reachable in exactly $k$ steps and the reachable values set $\rea$ defined as follows:
\begin{equation}
\label{reachseq}
\rea_k=\left\{\begin{array}{lr}
\displaystyle{A^{k}(\xin)+\sum_{i=0}^{k-1} A^i b} & \text{ if } k>0\\
\xin & \text{ if } k=0
\end{array}\right.
\quad\text{ and }\quad
\rea =\bigcup_{k\in\nn} \rea_k
\end{equation}
As $\xin$ is compact, the sets $\rea_k$ are compact for all $k\in\nn$. Finally, given a $d\times d$ symmetric matrix $Q$ and a $\rd$ vector $q$, we are interested in solving the following quadratic maximization problem :
\begin{equation}
\label{optpb}
 \sup_{x\in\rea} x^\intercal Q x +q^\intercal x=\sup_{k\in\nn} \max_{x\in \rea_k} x^\intercal Q x+q^\intercal x \enspace.
 \end{equation}
In the paper, we will use the following notations:
\begin{equation}
\label{optnuk}
\nu_k:=\max_{x\in \rea_k} x^\intercal Q x+q^\intercal x\quad \text{ and }\quad  
\nuopt:=\sup_{k\in\nn} \nu_k
\end{equation}
To compute $\nu_k$ is equivalent to solve a linearly constrained quadratic maximization problem. More or less efficient techniques are available to solve it as mentioned later in Subsection~\ref{quadprog}. The main issue is the fact that we have to solve potentially an infinite number of maximization problems. One way to compute exactly $\nuopt$ is to extract a finite number of them and prove that their resolution is sufficient to compute $\nuopt$. To extract this finite family of problem, we consider the following problem.
\begin{problem}
\label{problemtosolve}
Compute the smallest element of the set:
\begin{equation}
\label{eqpbinit}
\gse_\nu:=\{K\in\nn : \max_{k=0,\ldots,K} \nu_k\geq \sup_{j\geq K+1} \nu_j\}
\end{equation}
By convention, if the set $\gse_\nu$ is empty, its smallest element is $+\infty$.
\end{problem} 
To use an element $K$ in $\gse_\nu$ is interesting since it suffices to solve at most $K+1$ quadratic programs to compute $\nuopt$. The purpose of the paper is to construct a \emph{not too big} element of $\gse_\nu$ in order to solve a \emph{small number} of quadratic programs. If the constructed element $K$ is proved to be the smallest element of $\gse_\nu$, then $\nuopt$ can be computed with only one quadratic program for which the constraints set is $\rea_K$. The difficulty in the computation of an element of $\gse_\nu$ comes from its characterization using past and future values of $\nu_k$.
\subsection{Affine to Linear Systems}
\label{affinediscuss}
We introduced the problem where the dynamics is affine. Actually, it suffices to study the case where the dynamics is purely linear i.e. $b=0$ in Eq.~\eqref{stateeq}. In this case, the expression of $\nu_k$ becomes simpler:
\begin{equation}
\label{linearpb}
\nu_k=\max_{x\in \xin} f_k(x)\text{ where } x\mapsto f_k(x):=x^\intercal {A^\intercal}^k Q A^k x+q^\intercal A^k x
\end{equation}
We can adopt the basic approach which consists in using an auxiliary linear discrete-time system. If we suppose that $b\neq 0$ and $\Idd-A$ is invertible, then the sequence defined, for all $k\in\nn$, by $y_k=x_k-\btilde$ where $\btilde=(\Idd-A)^{-1} b$ satisfies, for all $k\in\nn$,  $y_{k+1}=A y_k$ and $x_k=A^k y_0+\btilde$. This latter expression leads to a new formulation of Problem~\eqref{optpb}: 
\[
\begin{array}{ll}
\displaystyle{\sup_{k\in\nn} \nu_k} &=\displaystyle{\sup_{k\in\nn} \max_{y_0\in\xin-\btilde} y_0^\intercal (A^k)^\intercal Q  A^k y_0+(2\btilde^\intercal Q+q^\intercal) A^k y_0 +\btilde^\intercal Q \btilde+q^\intercal \btilde}\\
   &=\left(\displaystyle{\sup_{k\in\nn} \max_{y_0\in\xin-\btilde} y_0^\intercal (A^k)^\intercal Q  A^k y_0+(2\btilde^\intercal Q+q^\intercal) A^k y_0}\right) +\btilde^\intercal Q \btilde+q^\intercal \btilde
    \end{array}
 \]
We recover an optimization problem of the form presented at Equation~\eqref{linearpb} up to an additive constant. The homogeneous of degree two part two is still $Q$ and the linear part becomes $2Q\btilde+q$. The polytope of initial conditions also changes since we have to consider now $\xin-\btilde$. 

In the rest of the paper, we will suppose that the system is purely linear and $b=0$ in Eq.~\eqref{stateeq}. When $k\in\nn$ is fixed, we have to compute $\nu_k$ defined at Eq.~\eqref{linearpb}. 
\subsection{Quadratic Programming}
\label{quadprog}
Once Problem~\ref{problemtosolve} solved, we have to solve a finite number of linearly constrained quadratic maximization problems. The complexity of the those maximization problems depends on the nature of the matrix $Q$. 

\paragraph{$Q$ negative semidefinite}
When $Q$ is negative semidefinite making the functions $x\mapsto x^\intercal {A^\intercal}^k Q A^k x+q^\intercal A^k x$ concave, for all $k\in\nn$, to compute $\nu_k$ can be done efficiently by any convex quadratic programming solver. Those solvers can be based on general non-linear methods, for example, interior points methods~\cite{vanderbei1999loqo,friedlander2012primal}, non-interior points methods~\cite{tian2015exterior} or active-sets methods~\cite{curtis2015globally,forsgren2016primal}.

\paragraph{$Q$ positive semidefinite} When $Q$ is positive semidefinite, the functions $x\mapsto x^\intercal {A^\intercal}^k Q A^k x+q^\intercal A^k x$ are convex, for all $k\in\nn$. It is well-known that there exists a vertex which maximizes the function. Then, as the number of vertices of $\xin$ is finite, we can find a maximizer and to compute $\nu_k$ with a finite number of evaluations by exploring all vertices of $\xin$.
Unfortunately, to use this method to solve a concave quadratic program does not scale well. Indeed, the number of vertices can increase exponentially as the dimension grows and the resolution time and the memory consumption blow up. More tractable methods to solve higher dimensional problems could be, for example, cutting plane methods~\cite{Tuy2016}; reductions to bilinear programming (with the use of particular cutting planes)~\cite{konno1976maximization} or branch-and-bound type methods~\cite{zamani2019new}. Some other methods can be found in the survey~\cite{Floudas1995}.

\paragraph{$Q$ indefinite case}
If $Q$ is indefinite, the functions $x\mapsto x^\intercal {A^\intercal}^k Q A^k x+q^\intercal A^k x$ are, in general, neither convex nor concave. Hence, for all $k\in\nn$, $\nu_k$ is the optimal value of an indefinite quadratic program. It is well-known (e.g see the survey in the paper~\cite{furini2019qplib}) that the problem is NP-hard. To provide an exact optimal solution to an indefinite quadratic program requires expensive methods such as branch-and-bound type algorithms~\cite{FLOUDAS19901397,burer2008finite}. More tractable algorithms and solvers can only guarantee local optimality of the returned solution~\cite{absil2007newton,huyer2018minq8}. Hopefully, in some specific situations, a global maximizer can be computed~\cite{zhang2000quadratic,kim2003exact,burer2019exact}.
\subsection{Previous Results from Diagonalizable case}
In~\cite{DBLP:journals/jota/Adje21}, we have proposed a method to compute an element of $\gse_\nu$ in the case where the matrix $A$ in Eq.~\eqref{stateeq} was diagonalizable. More precisely, the assumptions were:
\begin{enumerate}
\item $\rho(A)<1$;
\item There exists a non-singular complex matrix $U$ and a diagonal matrix $D$ such that $A=UDU^{-1}$;
\item There exists $k\in \nn$ such that $\nu_k>0$.
\end{enumerate}
With those assumptions, we have defined, for all $k\in\nn$ such that $\nu_k>0$:
\begin{equation}
\label{oldformula}
\kdiag_\nu(k):=\left\lfloor\dfrac{1}{\ln(\rho(A))}\ln\left(\dfrac{\sqrt{\nu_k+\vdiag^2}-\vdiag}{\sqrt{|\lmax{U^*QU}|\mus{(UU^{*})^{-1}}}}\right)\right\rfloor+1
\end{equation}
where
$\vdiag:=\dfrac{\norm{U^*q}_2}{2\sqrt{|\lmax{U^*QU}|}}$ and for all matrices Hermitian positive definite $B$, $\mu(B):= \displaystyle{\max_{x\in\xin} x^\intercal B x}$.

We have also got the following result:
\begin{theorem}[Th. 3.1 and Prop. 3.3 of ~\cite{DBLP:journals/jota/Adje21}] 
Let $j\in\nn$ such that $\nu_j>0$, we have:
\begin{enumerate}
\item $\max\{j,\kdiag_\nu(j)\}\in\gse_\nu$;
\item for all $k\geq \kdiag_\nu(j)$, $\nu_k\leq \nu_j$;
\item for all $k\in\nn$ such that $\nu_j\leq \nu_k$, $\kdiag_\nu(k)\leq \kdiag_\nu(j)$.
\end{enumerate}
\end{theorem}
We have provided an algorithm (Algorithm~\ref{algoAffQP}) to compute $\nuopt$, an integer $\kopt$ and a $\rd$-vector $\xopt$. The integer $\kopt$ represents a rank of a term $\nu_k$ that achieves the maximal term $\nuopt$. The vector $\xopt$ is an initial condition in $\xin$ that achieves the maximum in Eq.~\eqref{linearpb} for $k=\kopt$. However, the computation of $\nu_k$ requires an oracle capable of solving linearly constrained quadratic maximization problems. This oracle called SolveQP is used as ($\beta$,$x$)=SolveQP($g$,$\xin$) where $g$ is the quadratic objective function and $\xin$ is the polytope of initial conditions. Its outputs ($\beta$,$x$) are respectively the optimal value and an optimal solution. Additionally to the problem data, Algorithm~\ref{algoAffQP} requires an integer $N$ that stops the search of a positive term $\nu_k$. If no positive term is found before $N$ the method fails.
\begin{algorithm}[h!]
\DontPrintSemicolon

\Input{$Q$, $q$, $A$, $\xin$ and $N\in\nn$ a maximal integer for the search of a positive term.} 
\Output{A vector $(\nuopt,\kopt,\xopt)$ where $f_{\kopt}(\xopt)=\nuopt$ or a status "failed" if $\nu_k\leq 0$ for all $k=0,\ldots,N$.}
\Begin{
Compute an eigendecomposition $UDU^{-1}$ of $A$\;
Compute $\lmax{U^*QU}$, $\mu((UU^*)^{-1})$ and $\vdiag$\;
$(\nu_0,x_0)$=SolveQP($f_0$,$\xin$)\;
\If{$\nu_0\geq \left(\sqrt{|\lmax{U^*QU}|\mus{(UU^*)^{-1}}}+\vdiag\right)^2-\vdiag^2$}
{Return $(\nu_0,x_0,0)$} 
\Else
{
$k=0$\;
\While{$k<N$ and $\nu_k\leq 0$}
{
$k=k+1$\;
$(\nu_k,x_k)$=SolveQP($f_k$,$\xin$)\;
}

\If{$k=N$ and $\nu_k\leq 0$}{
Return "failed"
}
\Else{
$K=\kdiag_\nu(k)$\;
$\nuopt=\nu_k$; $\xopt=x_k$; $\kopt=k$\;
\While{$k<K$}{
$k=k+1$\;
$(\nu_k,x_k)$=SolveQP($f_k$,$\xin$)\;
\If{$\nuopt<\nu_k$}{
$\nuopt=\nu_k; \xopt=x_k; \kopt=k$\;
$K=\kdiag_\nu(k)$\;
}
}
Return $(\nuopt,\xopt,\kopt)$\;
}
}
}
\caption{Computation of the optimal value and an optimal solution for Eq.~\eqref{optpb}}
\label{algoAffQP}
\end{algorithm}

\begin{theorem}[Th. 3.5  and Prop. 3.6 of~\cite{DBLP:journals/jota/Adje21}]
\label{thold}
Suppose that there exists $k\leq N$ such that $\nu_k>0$. Let us denote $\Gamma(\nu):=\{k\in\nn : \nu_k>\nu_j,\ \forall j<k\}$. Then:
\begin{enumerate}
\item Then the sequence generated at Line 22 of Algorithm~\ref{algoAffQP}, $\left(\kdiag_\nu(k)\right)_{k\in \Gamma(\nu)}$ is decreasing;
\item Algorithm~\ref{algoAffQP} returns the optimal value $\nuopt$ and an optimal solution $(\xopt,\kopt)\in\xin\times \nn$ for $\sup\max$ problem depicted at Eq.~\eqref{optpb}. 
\end{enumerate}
\end{theorem}
Actually, the set $\Gamma(\nu)$ is finite and nonempty if $\gse_\nu$ is nonempty and we have $\max \Gamma(\nu) = \min \gse_\nu$. The first statement of Th.~\ref{thold} proves that we reduce the number of iterations by computing a new $\kdiag_\nu(k)$ at Line 22 of Algorithm~\ref{algoAffQP}. The second statement proves that Algorithm~\ref{algoAffQP} is correct.

The main goal of this paper is to generalize the approach to a  nondiagonalizable matrix $A$. The assumption $\rho(A)<1$ and the existence of a positive term are still valid in this paper. We will keep the same algorithm and use a formula similar to the one depicted at Eq.~\eqref{oldformula}. 
\subsection{Running Example}
\label{runningex}
To illustrate the results, we will consider the following family of linear systems for which the system matrix is known to be not diagonalizable:
\begin{equation}
\label{runningdata}
\xin=[-1,1]^2 \quad \text{ and } \quad
A_\gamma=\gamma\begin{pmatrix} 1 & 1 \\ 0 & 1\end{pmatrix},\ 
\end{equation}
where $\gamma\in (0,1)$ in order to satisfy the stability assumption $\rho(A_\gamma)=\gamma<1$. Moreover, we will illustrate the optimization problems with a very simple homogeneous quadratic objective function characterized by:
\begin{equation}
\label{runningobj}
\qex=\begin{pmatrix} 1 & 0 \\ 0 & 0\end{pmatrix}
\end{equation}
As $\qex$ is semidefinite positive, the functions $x\mapsto x^\intercal {A^k}^\intercal Q A^k x$ are convex for all $k\in\nn$ and to maximize it on $\xin$, it suffices to consider the vertices of $\xin$. For all $k\in\nn$, we have :
\[
A_\gamma^k=\gamma^{k}\begin{pmatrix}
1 & k \\ 0 & 1
\end{pmatrix}
\]
and
\[
\nu_k=\max_{x \in [-1,1]^2} x^\intercal {A_\gamma^k}^\intercal \qex A_\gamma^k x=\gamma^{2k}(1+k)^2
\]
We conclude that $\nuopt=\max_{k\in\nn} \gamma^{2k}(1+k)^2$. Let us introduce the function $\psi_\gamma:u\mapsto \gamma^{2u}(1+u)^2$ and
\[
k_1^*=\left\lfloor \dfrac{1}{-\ln(\gamma)}-1\right\rfloor\text{ and } 
k_2^*=\left\lfloor \dfrac{1}{-\ln(\gamma)}\right\rfloor
\]
Let us suppose that $\gamma> \exp(-1)$. In this case, the integers $k_1^*$ and $k_2^*$ surrond the maximizer of $\psi_\gamma$ over the reals which is $u^*:=-1/\ln(\gamma)-1$. As $\psi_\gamma$ increases on $(0,u^*)$ and decreases on $(u^*,+\infty)$, the integers $k_1^*$ and $k_2^*$ are the only possible integral maximizers of $\psi_\gamma$. Hence:
\[
\left\{
\begin{array}{ll}
\nuopt=1, \kopt=0& \text{if } \gamma\leq \exp(-1)\\
(\nuopt,\kopt)=\left\{
\begin{array}{ll}
(\psi_\gamma(k_1^*),k_1^*) & \text{ if } \psi_\gamma(k_1^*)\geq \psi_\gamma(k_2^*) \\
(\psi_\gamma(k_2^*),k_2^*) & \text{ if } \psi_\gamma(k_2^*)>\psi_\gamma(k_1^*) \\
\end{array} \right. & \text{ if }\gamma > \exp(-1)
\end{array}
\right.
\]
We remark that, when $\gamma >\exp(-1)$, $\kopt$ and $\nuopt$ are increasing as functions of $\gamma$. Moreover, they tend to $+\infty$ as $\gamma$ tends to 1.
As $(1-x)^{-1}-1< (-\ln(x))^{-1}<(1-x)^{-1}$ for all $x\in (0,1)$, we conclude that :
\[
k_1^*\in\left\{\left\lfloor\dfrac{1}{1-\gamma}-2\right\rfloor,\left\lfloor\dfrac{1}{1-\gamma}-1\right\rfloor\right\} \text{ and } k_2^*\in\left\{\left\lfloor\dfrac{1}{1-\gamma}-1\right\rfloor,\left\lfloor \dfrac{1}{1-\gamma}\right\rfloor\right\} 
\]
Now to obtain the maximizer on $\nn$, we have to compare $\psi_\gamma(k_1^*)$ and
$\psi_\gamma(k_2^*)$ or more generally $ \psi_\gamma(u+1)$ and $\psi_\gamma(u)$. The inequality $\psi_\gamma(u+1)/\psi_\gamma(u)\geq 1$ is equivalent to $(1-\gamma)(u+1)\leq \gamma$ which is the same as $u\leq \gamma (1-\gamma)^{-1}-1=(\gamma-1)(1-\gamma)^{-1}+(1-\gamma)^{-1}-1=(1-\gamma)^{-1}-2$. From the double inequalities  $(1-x)^{-1}-1< (-\ln(x))^{-1}< (1-x)^{-1}$ for all $x\in (0,1)$, we have $k_1^*< (1-\gamma)^{-1}-2$ if and only if $k_1^*=\lfloor(1-\gamma)^{-1}-2\rfloor$ and $(1-\gamma)^{-1}-2$ is not an integer. In this case, $k_2^*$ is the optimal solution. The equality 
$\psi_\gamma(k_1^*)=\psi_\gamma(k_2^*)$ holds if and only if $k_1^*=(1-\gamma)^{-1}-2$ (which implies that $(1-\gamma)^{-1}-2)$ is an integer). Finally, $k_1^*$ is the only optimal solution when $k_1^*=\lfloor (1-\gamma)^{-1}\rfloor-1$. In summary:
\[
\kopt=\left\{ 
\begin{array}{lr}
k_1^* & \text{if } k_1^*=\lfloor (1-\gamma)^{-1}\rfloor-1\\
k_2^* & \text{if } k_1^*=\lfloor (1-\gamma)^{-1}\rfloor-2 \text{ and } (1-\gamma)^{-1}\notin \nn
\\
\{k_1^*,k_2^*\} & \text{if } k_1^*=(1-\gamma)^{-1}-2
\end{array}
\right.
\]
To illustrate our words with numerical experiments, we will consider the four following families :
\begin{itemize}
\item The family $\gamma_n:=\exp(-1/n)$ parameterized by $n\in\nn^*$ satisfies 
$\kopt=k_1^*=n-1$ and $\nuopt=\gamma_n^{2(n-1)}n^2$. Indeed, $k_1^*=\lfloor (-\ln(\exp(-1/n)))^{-1} -1 \rfloor=(-\ln(\exp(-1/n)))^{-1} -1=n-1$. We conclude that $k_1^*$ is also equal to $\lfloor (1-\exp(-1/n))^{-1}\rfloor-1$ by the fact that for $x\in\rr$ and $z\in\nn$ such that $x-1<z<x$, we must have $z=\lfloor x\rfloor$. 
\item Let us take the family, $\gamma_n:=(n+0.1)/(n+1+0.1)$ where $n\in\nn^*$. 
For this family, we need a tighter inequality on logarithm~\cite{b0750f61-9854-3f22-8cec-fb4c23a27312}, for all $x>0$:
$x^{-1}\leq (\ln(1+x))^{-1}\leq (2+x)(2x)^{-1}$. We need this as, in this case, $-\ln(\gamma_n)=\ln(n+1.1/n+0.1)=\ln(1+(n+0.1)^{-1})$ and so $n-1=\lfloor n+0.1\rfloor -1\leq k_1^*\leq \lfloor n+0.6\rfloor-1=n-1$. We conclude that, for this family, $k_2^*$ is the only optimal solution since $\lfloor (1-\gamma_n)^{-1}\rfloor-2=\lfloor n+1.1\rfloor -2 =n-1=k_1^*$ and $(1-\gamma_n)^{-1}=n+1.1\notin\nn$.
\item Let us consider the family $\gamma_n:=1-10^{-n}$ indexed by $n\in\nn^*$. From the inequalities $(1-x)^{-1}-1< -\ln(x)^{-1}< (1-x)^{-1}$, for all $x\in (0,1)$, we obtain $k_1^*=10^n-2=(1-\gamma_n)^{-1}-2$ and $k_2^*=10^n-1$. Indeed, when $y-1<z<y$ for $y\in\nn$ and $z\in\rr$, we must have $\lfloor z\rfloor=y-1$. 
We conclude that, for this family, $k_1^*$ and $k_2^*$ are both optimal.  
\item Finally we consider the family, $\gamma_n:=1-10^{-n}+10^{-n-1}$ indexed by $n\in\nn^*$. We use a coarser inequality $(1-x)^{-1}-1< -\ln(x)^{-1}< (2(1-x))^{-1}(1+x)$ for all $x>1$. Hence, we have $k_2^*=\lfloor (10^{-n}-10^{-n-1})^{-1}-1 \rfloor=\lfloor 10^{n}/(1-10^{-1})\rfloor-1=\lfloor \lim_{M\to \infty}10^{n}(1-10^{-(M-n+1)})/(1-10^{-1})\rfloor-1=\lfloor \lim_{M\to -\infty} \sum_{i=M}^n 10^{i}\rfloor-1=\sum_{i=1}^n 10^{i}$. We conclude that $k_1^*=\lfloor (1-\gamma_n)^{-1}\rfloor-2$, $(1-\gamma_n)^{-1}\notin\nn$ and $k_2^*$ is the only optimal solution.
\end{itemize}
\subsection{On the Supremum of Zero Limit Real Sequences}
\label{sequences}
Let us denote by $\co$ the set of real sequences the limit of which is equal to 0 i.e. $\co=\{u=(u_0,u_1,\ldots)\in\rr^\nn: \lim_{n\to +\infty} u_n=0\}$. For an element of $\co$, we are interested in computing the supremum of its terms and a maximizer as well. We will write $\agmx(u):=\{k\in\nn: u_k=\sup_{n\in\nn} u_n\}$. For $u\in\co$, we can characterize $\agmx(u)$ from the two sets of ranks:
\[
\gse_u:=\left\{k\in\nn: \max_{0\leq j\leq k} u_j\geq \sup_{j>k} u_j\right\}\ \text{ and }\ \gs_u:=\left\{k\in\nn: \max_{0\leq j\leq k} u_j>\sup_{j>k} u_j\right\}
\]
It is obvious that $\gse_u\subseteq \gs_u$. Moreover, if $k$ belongs to $\gse_u$ (resp. $\gs_u$), then any integer greater than $k$ belongs to $\gse_u$ (resp. $\gs_u$). Besides, the emptiness of $\agmx(u)$ depends on the existence of non-negative and positive terms:
\[
\poseq_u=\{k\in\nn: u_k\geq 0\}\qquad \text{ and }\qquad \pos_u=\{k\in\nn: u_k>0\}\enspace .
\]
We have readily $\pos_u\subseteq \poseq_u$. Note that even if $u\in\co$, $\poseq_u$ and $\pos_u$ can be empty. Finally, we need the smallest elements of the latter sets of ranks.
\[
\bigkse_u=\inf \gse_u ;\  
\qquad 
\bigks_u=\inf\gs_u ;\
\qquad
\kse_u=\inf\poseq_u
\quad \text{ and }\ 
\ks_u=\inf\pos_u
\]
Note that by convention, the smallest element of the empty set is equal to $+\infty$.
\begin{proposition}
\label{posetdelta}
Let $u\in\co$. The following assertions hold:
\begin{enumerate}
    \item For all $k\in\nn$, $\sup_{l\geq k} u_l\geq 0$;
    \item $\poseq_u\neq \emptyset\iff \gse_u\neq\emptyset$;
    \item $\pos_u\neq \emptyset\iff \gs_u\neq\emptyset$;
    \item $\gs_u=\emptyset\iff \sup_{k\in\nn} u_k=0$.
        \end{enumerate}
\end{proposition}
\begin{proposition}
\label{ineg}
The following inequalities hold:
\[
\kse_u\leq \ks_u;\qquad \kse_u\leq \bigkse_u;\qquad \ks_u\leq \bigks_u\quad \text{ and }\bigkse_u\leq \bigks_u\enspace .
\]
Moreover, if $\ks_u<+\infty$, then $\ks_u\leq \bigkse_u$.
\end{proposition}
\begin{proposition}[Argmax]
\label{argmax}
Let $u\in\co$. The following assertions hold:
\begin{enumerate}
    \item $\agmx(u)\subseteq \gse_u$;
    \item If $\bigkse_u<+\infty$, $\bigkse_u=\min \agmx(u)$;
    \item If $\bigks_u<+\infty$ then $\sup_{k\in\nn} u_k>0$ and $\bigks_u=\max \agmx(u)$;
    \item $\agmx(u)\neq\emptyset\iff \gse_u\neq\emptyset$;
    \item $\agmx(u)=\emptyset \implies \sup_{n\in\nn} u_n=\lim_{k\to +\infty} u_k=0$; 
\end{enumerate}
\end{proposition}
\begin{proposition}
\label{supremumpos}
Let $u\in\co$. The following assertions hold:
\begin{enumerate}
    \item Assume that $\poseq_u\neq \emptyset$. Let $k\geq \kse_u$ such that $\displaystyle{\max_{\kse_u\leq j\leq k} u_j\geq \sup_{j>k} u_j}$ then $k\in\gse_u$;
    \item Assume that $\pos_u\neq \emptyset$. For all $\bigkse_u\leq k$, 
\[    
    \sup_{0\leq j\leq k} u_j=\sup_{\ks_u\leq j\leq k} u_j=\sup_{\kse_u\leq j\leq k} u_j=\sup_{j\in\nn} u_j\enspace. 
    \]
\end{enumerate}
\end{proposition}
In summary, for $u\in\co$, to compute $\sup_{k\in\nn} u_k$, we need to study first the emptiness of $\pos_u$. If $\pos_u$ is empty, we know that (the fourth statement of Prop.~\ref{posetdelta}) $\sup_{k\in\nn} u_k$ is equal to 0. If $\pos_u$ is not empty, then we have to compute $\bigkse_u$ to find $\sup_{k\in\nn} u_k$ and a maximizer (the second statement of Prop.~\ref{argmax}). The decision problem about the existence of a positive term for a sequence is still open~\cite{ouaknine2014positivity}.  Moreover, to identify $\bigkse_u$, we need knowledge on the past and the future of the sequence. The good point is that any overapproximation $k$ of $\bigkse_u$ permits to know $\sup_{k\in\nn} u_k$ by computing $\sup_{0\leq j\leq k} u_j$ (the third statement of Prop~\ref{supremumpos}).
\section{The Key Formula To Overapproximate A Rank Maximizer}
\label{mainresults}
We come back to the $\sup\max$ problem described at Eq.~\eqref{optpb}. We recall that from the discussion of Subsection~\ref{affinediscuss}, we suppose that $b=0$. Moreover, with $b=0$, if $\xin=\{0\}$, we have obviously $\sup_{k\in\nn}\nu_k=0$. We also remark that if $Q\preceq 0$ and $q=0$, from Prop.~\ref{posetdelta}, the problem is trivially solved with $\nuopt=0$. Finally, $\rho(A)=0$ if and only if $A$ is nilpotent. In this case, we must have $A^d=0$, we conclude that $\sup_{k\in\nn} \nu_k=\sup_{0\leq k\leq d} \nu_k$. In this case, the problem is solved and we suppose that $\rho(A)>0$. Note that if $A$ is nilpotent and not null, it is not diagonalizable and this assumption was not necessary in~\cite{DBLP:journals/jota/Adje21}. 

Let us consider the initial problem :
\begin{equation}
\label{eqnuk}
\nu_k=\sup_{x\in\xin} x^\intercal {A^k}^\intercal Q A^k x+q^\intercal A^k x
\end{equation}
In the rest of the paper, we assume, to avoid trivially solved problems:
\begin{enumerate}
\item If $Q\preceq 0$, then $q\neq 0$;
\item $\xin\neq \{0\}$;
\item $\rho(A)>0$.
\end{enumerate}
We also make the following assumption:
\begin{assumption}
\label{assum1}
The matrix $A$ is suppose to be stable i.e. $\rho(A)<1$. With the previous assumption, we suppose that the spectral radius of $A$ satisfies $0<\rho(A)<1$. 
\end{assumption}
Recall that $\xin$ is a polytope and then is bounded, the following proposition thus holds.
\begin{proposition}
\label{propnuk0}
Assumption~\ref{assum1} ensures that $\nu\in\co$. 
\end{proposition}

\subsection{Tools from Linear Algebra}
First, we recall that, for a Hermitian matrix $B\in\cdd$ (i.e. $B$ satisfies $B^*=B$), $x^* B x$ is real for all $x\in\cd$; the eigenvalues and the diagonal elements are also real. The set of Hermitian matrices will be denoted  by $\hdd$. We use the standard notation $B\succ 0$ (resp. $B\succeq 0$) for $B$ being Hermitian positive definite (resp. Hermitian positive semidefinite) meaning that $x^* B x>0$ for all non-zero $x\in\cd$ (resp. $x^* B x\geq 0$ for all $x\in\cd$). We also use the notation $\hddp$ (resp. $\hp$) for 
the set of positive definite (resp. positive semidefinite) matrices.

Let us introduce the set $\lyap{A}$ of solutions of the discrete Lyapunov equation with respect to $A$ i.e.
\[
\lyap{A}=\{P\in\hdd: P\succ 0 \text{ and } P-A^\intercal P A\succ 0\}\enspace .
\]
Assumption~\ref{assum1} ensures that $\lyap{A}\neq \emptyset$. Moreover,  elements of $\lyap{A}$ can be characterized by closed formula: 
\[
P\in\lyap{A}\iff \text{ there exists } Q\succ 0\text{ s.t. } P=\sum_{j=0}^{\infty} (A^j)^\intercal Q A^j
\]
\begin{example}
\label{exampleone}
Let us consider the family of matrices $\{A_\gamma\}_{\gamma\in(0,1)}$ introduced at Eq.~\eqref{runningdata}.  Let $\gamma\in (0,1)$. To simplify the computations, we will interest in diagonal matrices in $\lyap{A_\gamma}$. Then we consider a closed strict subset of $\lyap{A_\gamma}$. We write, for $a,b\in\rr$, the matrix $\pab=\begin{pmatrix} a & 0 \\ 0 & b\end{pmatrix}$. Then $\pab$ belongs to $\lyap{A_\gamma}$ iff the conditions are fulfilled:
\begin{enumerate}
\item $\pab\succ 0\iff a > 0 \text{ and } b>0$;
\item $\pab-A_\gamma^\intercal \pab A_\gamma\succ 0$.
\end{enumerate}
Those two latter conditions are, respectively, equivalent to 
$a$ and $b$ are strictly positive and to $a$ is strictly positive and $a(1-\gamma^2)(b(1-\gamma^2)-\gamma^2 a)-\gamma^4 a^2$ is strictly positive. Finally :
\begin{equation}
\label{lyaprunn}
\pab\in\lyap{A_\gamma}\iff \left\{ \begin{array}{l}a>0\\ b>\dfrac{\gamma^2}{(1-\gamma^2)^2} a\end{array}\right. 
\end{equation}
We introduce the subset of  $\rr^2$: \[\lyapex(\gamma):=\{(x,y)\in\rr^2 : x>0,\ y>\gamma^2(1-\gamma^2)^{-2}x\}\] By doing so, $\pab\in\lyap{A_\gamma}\iff (a,b)\in\lyapex_\gamma$.
\qed
\end{example}
We recall that a matrix norm is real-valued function $N$ on $\cdd$ that satisfies for all $B,C\in\cdd$ and all $\alpha\in\cc$ the three following statements: 
\begin{inparaenum}[(i)]
\item $N(B)\geq 0$ and $N(B)=0\iff B=0$; 
\item $N(\alpha B)=|\alpha| N(B)$; 
\item $N(B+C)\leq N(B)+N(C)$ 
\end{inparaenum}.
If in addition :
\begin{inparaenum}[(i)]
\setcounter{enumi}{3}
\item $N(BC)\leq N(B) N(C)$
\end{inparaenum}
$N$ is said to be a submultiplicative matrix norm. A submultiplicative matrix norm $N$ verifies $N(B^k)\leq N(B)^k$ for all $B\in\cdd$ and all $k\in\nn^*$. All matrix norms on $\cdd$ are equivalent i.e. for two matrix norms $N,N'$, there exists two strictly positive scalars $\alpha$, $\beta$ s.t. $\alpha N \leq N'\leq \beta N$ (in the functional order sense). We recall that $B\in\hddp$ has a (unique) square root $B^{1/2}$ defined by $B^{1/2}=UD^{1/2}U^*$ where $B=UDU^*$ and $UU^* =U^*U=\Idd_d$. We also have $B^{-1/2}=(B^{-1})^{1/2}=(B^{1/2})^{-1}$.
 
Since for all $P\in\lyap{A}$, $P\succ 0$, we can define the weighted norm $\norm{\cdot}_P$ on $\cd$ defined for all $x\in\cd$ by $\norm{x}_P:=\sqrt{x^* P x}$. From it, we can define the associated operator matrix norm of $A$ :
\begin{equation}
\label{normdef}
\begin{array}{lcl}
    \norm{A}_P^2&=&\displaystyle{\max_{x\neq 0} \dfrac{x^* A^* P A x}{x^* P x}}\\
    &=&\displaystyle{\max_{x^* P x=1} x^* A^* P A x}\\
    &=&\displaystyle{\inf\{\alpha>0\mid \alpha P-A^* P A\succeq 0\}}\\
    &=&\displaystyle{\lmax{P^{-1/2}A^*P AP^{-1/2}}}\\
    &=&\displaystyle{\lmax{A^*PAP^{-1}}}
    \end{array}
\end{equation}
Operator matrix norms are submultiplicative and their value at $\Idd$ equals to one. Hence, for all $P\in\lyap{A}$ and for all $k\in\nn$, $\norm{A^k}_P\leq \norm{A}_P^k$. From~\eqref{normdef}, $\norm{A}_P$ can be computed using either semidefinite programming and associated solvers or classical numerical linear algebra tools to compute the maximal eigenvalue of a complex matrix, Hermitian matrix or not.
\begin{lemma}
\label{lyapnorm}
The following statement are valid:
\begin{enumerate}
\item The map defined on $\hdd$, $P\mapsto P-A^\intercal P A$ is linear and thus continuous;
\item For all $\alpha>0$, $\norm{A}_{\alpha P}=\norm{A}_P$;
\item For all $P\in\lyap{A}$, $\rho(A)\leq \norm{A}_P<1$.
\end{enumerate}
\end{lemma}
We also introduce the dual norm of $\norm{\cdot}_P$ denoted by $\norm{\cdot}_P^*$ and defined for all $y\in\cd$ by:
\[
\norm{y}_P^*=\max_{\norm{z}_P=1} |y^* z|=\norm{P^{-1/2}y}_2=\sqrt{y^* P^{-1} y} \enspace .
\]
We thus have for all $y,z\in\cd$, $|y^* z|\leq \norm{y}^{*}_P\norm{z}_P$.
\begin{example}
\label{examplenorm}
Let $(a,b)\in\lyapex(\gamma)$. We use the formula $\lmax{A^* \pab   A\pab^{-1}}$ to compute $\norm{A_\gamma}_{\pab}$. We have simply $\pab^{-1}=\begin{psmallmatrix} a^{-1} & 0\\ 0 & b^{-1}\end{psmallmatrix}$. We obtain:
\[
\norm{A_\gamma}_{\pab}^2=\dfrac{\gamma^2}{2}\left(2+\dfrac{a}{b}+\sqrt{\dfrac{a}{b}\left(\dfrac{a}{b}+4\right)}\right)
=\dfrac{\gamma^2}{2}\left(\dfrac{1+2\frac{b}{a}+\sqrt{1+4\frac{b}{a}}}{\frac{b}{a}}\right)
\]
Applying the change of variable $b'=b/a>\gamma^{2}(1-\gamma^2)^{-2}$, we have $(1,b')\in\lyapex(\gamma)$ and we conclude that:
\[
\norm{A_\gamma}_{\pab}^2=\norm{A_\gamma}_{\pubp}^2=\dfrac{\gamma^2}{2}\left(\dfrac{1+2b'+\sqrt{1+4b'}}{b'}\right)
\]
First, the function $x\mapsto x^{-1}(1+2x+\sqrt{1+4x})$ is continuous and strictly decreases. Second, its value at $x=\gamma^{2}(1-\gamma^2)^{-2}$ is equal to $2\gamma^{-2}$. Hence $\norm{A_\gamma}_{\pubp}^2$ is strictly smaller than 1. Finally the limit of $x\mapsto x^{-1}(1+2x+\sqrt{1+4x})$ at $+\infty$ is 2 and we conclude that $\norm{A_\gamma}_{\pubp}^2$ is greater than $\gamma^2=\rho(A_\gamma)^2$.
\qed
\end{example}
An element of $\lyap{A}$ is not necessarily greater than $Q$. Then, we will need to scale an element of $\lyap{A}$ in order to dominate $Q$. Hence, for a Hermitian matrix $P\in\hdd$, we introduce the following set:
\begin{equation}
\label{contPeq}
\cont(P)=\left\{t\in\rr_+ : tP-Q\succeq 0\right\}
\end{equation}
The next lemma aims to identify the smallest element of  $\cont(P)$. 
\begin{lemma}
\label{propcont}
Let $P$ be a Hermitian matrix. We define:
\begin{equation}
\label{deftopt}
\topt{P}:=\max\left\{\lmax{P^{-1/2} Q P^{-1/2}},0\right\}=\max\{\lmax{QP^{-1},0}\}
\end{equation}
The following statements hold :
\begin{enumerate}
\item $Q\preceq 0\iff 0\in \cont(P)$ ;
\item If $Q\npreceq 0$, $\cont(P)\subseteq \rr_+^*$;
\item For all $P\succ 0$, $\Min\cont(P)=\topt{P}$.
\end{enumerate}
\end{lemma}
\begin{example}
\label{examplelmaxt}
Recall that $\qex=\begin{psmallmatrix} 1 & 0\\ 0 & 0\end{psmallmatrix}$ and assume that $(a,b)\in\lyapex(\gamma)$. Since $Q\succeq 0$, 
$t_{\qex}^*(\pab)=\lmax{\qex\pab^{-1}}=a^{-1}$.
Finally for all $(a,b)\in\lyapex(\gamma)$, $t_{\qex}^*(\pab)\pab=\begin{psmallmatrix}
1 & 0\\ 0 & b a^{-1}\end{psmallmatrix}\succeq \qex$ and we cannot expect a smallest scaling for $\pab$.
\qed
\end{example}
Finally, we introduce the last set:
\begin{equation}
\label{ptset}
\ptset:=\left\{(P,t)\in\cdd\times \rr_+ : P\in\lyap{A} \text{ and }t\in\cont(P)\right\}
\end{equation}
\subsection{The Rank Formula for the General Case}
With those basic tools, we can provide an upper-bound over the terms of the sequence $(\nu_k)_{k\in\nn}$. Recall that, for all $B\succ 0$:
\[
\mu(B)=\max_{x\in\xin} x^\intercal B x
\]
\begin{example}
\label{examplemup}
Let $(a,b)\in\lyapex(\gamma)$. For $\xin=[-1,1]^2$, we have $\mu(\pab)=a+b$.
\qed
\end{example}
Let us introduce, for $(P,t)\in\ptset$, the function on $\rr_+$:
\begin{equation}
\label{homeoseq}
x\mapsto H_{P,t}(x)=\left\{
\begin{array}{lr}
x\sqrt{\mu(P)}\norm{q}_P^*  & \text{ if } t=0\\
\left(\sqrt{t\mu(P)}x+\dfrac{\norm{q}_P^*}{2\sqrt{t}}\right)^2-\dfrac{{\norm{q}_P^*}^2}{4t} & \text{ otherwise}
\end{array}
\right.
\end{equation}
The values $H_{P,t}(x)$ are well defined and not null for all strictly positive $x$. Indeed, the value $\mu(P)$ cannot be null as $P\in\lyap{A}$ and $\xin$ is not reduced to $\{0\}$.  Following Lemma~\ref{propcont}, when $Q\npreceq 0$, $0$ does not belong to $\cont(P)$ and $t\neq 0$ in Eq.~\eqref{homeoseq}. This also implies that the case $t=0$ can only be used when $Q\preceq 0$.
Finally since if $Q\preceq 0$, by assumption, we have $q\neq 0$, $\norm{q}_P^*$ cannot be null when $t=0$. 
\begin{proposition}
\label{homeoprop}
Let $(P,t)\in\ptset$. The following assertions are true:
\begin{itemize}
\item $H_{P,t}(x)=0$ if and only if $x=0$;
\item $H_{P,t}$ is strictly increasing;
\item $H_{P,t}$ is an homeomorphism from $\rr_+$ to itself and 
\[x\mapsto H^{-1}_{P,t}(x)=\left\{
\begin{array}{lr}
x\left(\sqrt{\mu(P)}\norm{q}_P^*\right)^{-1}  & \text{ if } t=0\\
\dfrac{1}{2t\sqrt{\mu(P)}}\left(\sqrt{4t x+(\norm{q}^*_P)^2}-\norm{q}^*_P \right) & \text{ otherwise}
\end{array}
\right.
\]
\item for all $\alpha>0$, $H_{P,t}=H_{\alpha P,\alpha^{-1}t}$ and $H^{-1}_{P,t}=H^{-1}_{\alpha P,\alpha^{-1} t}$.
\end{itemize}
\end{proposition}
This is simple to see that $H_{P,t}^{-1}$ is also strictly increasing as the inverse of a strictly increasing function.
\begin{proposition}
\label{nukineg}
Let $(P,t)\in \ptset$. We have :
\[
\nu_0\leq H_{P,t}(1)\text{ and for all } k\in\nn^*,\ \nu_k\leq H_{P,t}(\norm{A}_p^k)<H_{P,t}(1).
\]
\end{proposition}
\begin{proof}
Let $(P,t)\in \ptset$ and $x\in\xin$. Let suppose that $k=0$. If $t=0$, we have $Q\preceq 0$ and so $x^\intercal Qx+q^\intercal x\leq \norm{q}_P^* \norm{x}_P\leq \norm{q}_P^*\sqrt{\mu(P)}=H_{P,0}(1)$. Now, if $t>0$, by definition of $\ptset$, $\norm{\cdot}_P^*$ and $\mu(P)$, $x^\intercal Q x+q^\intercal x\leq tx^\intercal P x +\norm{q}_P^* \norm{x}_P\leq t\mu(P)+\norm{q}_P^*\sqrt{\mu(P)}=H_{P,t}(1)$. In both cases, taking the max over $\xin$ leads to $\nu_0\leq H_{P,t}(1)$.

Assume that $k\in\nn^*$. Then, we have $q^\intercal {A^k}^\intercal x\leq \norm{q}_P^* \norm{A^k x}_P\leq \norm{q}_P^* \norm{A^k}_P\norm{x}_P
\leq \norm{q}_P^* \norm{A}^k_P\sqrt{\mu(P)}$. If $t=0$ then $Q\preceq 0$ and $x^\intercal {A^k}^\intercal Q A^k x\leq 0$. Hence, $\nu_k\leq \max_{x\in\xin} q^\intercal A^k x\leq \norm{q}_P^* \norm{A}^k_P\sqrt{\mu(P)}=H_{P,0}(\norm{A}_P^k)$. Now suppose that $t>0$. We have
$x^\intercal {A^k}^\intercal Q A^k x\leq t\norm{A^k x}_P^2
\leq t\norm{A^k}_P^2 \mu(P)\leq t\norm{A}^{2k}_P \mu(P)$. By summation, we get :
\[
\nu_k=\max_{x\in \xin} x^\intercal {A^k}^\intercal Q A^k x+q^\intercal {A^k}^\intercal x\leq 
\left(\sqrt{t\mu(P)}\norm{A}_P^k+\dfrac{\norm{q}_P^*}{2\sqrt{t}}\right)^2-\dfrac{{\norm{q}_P^*}^2}{4t}
\enspace .
\]
The second strict inequality follows from Lemma~\ref{lyapnorm} and the fact that $H_{P,t}$ is strictly increasing.
\end{proof}
\begin{corollary}
\label{corollary1}
If there exists $(P,t)\in \ptset$ such that $\nu_0= H_{P,t}(1)$ then $\nu_0=\displaystyle{\sup_{k\in\nn} \nu_k}$. 
\end{corollary}
Corollary~\ref{corollary1} holds for example when $Q\in\lyap{A}$ and $q=0$. In this case, $(Q,1)\in\ptset$ and $\nu_0=\mu(P)$ with $P=Q$. Now, for the rest of the paper, we suppose that the hypothesis of Corollary~\ref{corollary1} does not hold i.e. :
\begin{equation}
\label{ineqforall}
\tag{A}
\forall\, (P,t)\in\ptset;\ \nu_0<H_{P,t}(1)
\end{equation}
Assumption~\eqref{ineqforall} implies that the strict inequality $\nu_k<H_{P,t}(1)$ is valid for all $(P,t)\in\ptset$ and all $k\in\nn$. In the rest of the paper, we assume that Assumption~\eqref{ineqforall} holds.
\begin{example}
\label{exampleprop8}
We illustrate Prop~\ref{nukineg} for $k=0$ and Assumption~\eqref{ineqforall} for our running example. We have $\nu_0=\max_{x\in [-1,1]^2} x^\intercal \qex x=1$. Let $(a,b)\in\lyapex(\gamma)$ and $t\in\cont(\pab)$. Suppose that $\nu_0=H_{\pab,t}(1)=t\mu(\pab)$. This is equivalent to $1=t(a+b)=ta(1+b/a)$. Since $t\in\cont(\pab)$, by Lemma~\ref{propcont}, we have $t\geq t_{\qex}^*( \pab)=a^{-1}$ and $ta\geq 1$. As $a$ and $b$ are strictly positive, $t(a+b)\geq 1+b/a$. We conclude that $1=t(a+b)\geq 1+b/a$ implies that $b/a\leq 0$ which contradicts the fact that $a$ and $b$ are strictly positive and $\nu_0<H_{\pab,t}(1)$.
\end{example}
\begin{remark}
Prop.~\ref{nukineg} includes the case where the objective function of the maximization problem depicted at Eq.~\eqref{optpb} is linear i.e. $Q=0$. 
\end{remark}
\begin{example}
\label{complicated}
Let $(a,b)\in\lyapex(\gamma)$ and $t\in\cont(\pab)$. Prop.~\ref{nukineg} applied to $(\pab,t)$ and $\nu_k=\gamma^{2k}(1+k)^2$ reads :
\begin{equation}
\label{eqexineq}
\gamma^{2k}(1+k)^2\leq t(a+b)\left(\dfrac{\gamma^2}{2}\left(\dfrac{a+2b+\sqrt{a(a+4b)}}{b}\right)\right)^{k}
\end{equation}
for all $k\in\nn$.

Since, writing as in Example~\ref{examplenorm}, $b'=b/a$ and $t(a+b)\geq 1+b'$ (as shown in Example~\ref{exampleprop8}):
\[
\begin{array}{ll}
t(a+b)\left(\gamma^2\dfrac{a+2b+\sqrt{a(a+4b)}}{2b}\right)^{k}&=ta(1+b')\left(\gamma^2\dfrac{(1+2b'+\sqrt{1+4b'})}{2b'}\right)^{k}\\
&\geq (1+b')\left(\gamma^2\dfrac{(1+2b'+\sqrt{1+4b'})}{2b'}\right)^{k}
\end{array}
\]
and $(P_{1,b'},1)\in \ptex$, Eq.~\eqref{eqexineq} is equivalent to:
\[
(1+k)^2\leq (1+b')\left(\dfrac{1+2b'+\sqrt{1+4b'}}{2b'}\right)^{k}
\]
We introduce $B:=(1+2b'+\sqrt{1+4b'})/(2b')$ and thus $b'=B/(B-1)^2$. As $B=\gamma^{-2}\norm{A_\gamma}_{P_{1,b'}}$, then $B\in (1,\gamma^{-2})$ and Eq.~\eqref{eqexineq} is equivalent to:
\[
(1+k)^2\leq \left(1+\dfrac{B}{(B-1)^2}\right)B^{k}
\]
for all $k\in\nn$, for all $B\in (1,\gamma^{-2})$. The latter inequality is not trivial and we propose an independent proof for it in appendix.

At Prop.~\ref{nukineg}, we also have $\nu_k<H_{P,t}(1)$ for all $k\in\nn$ (since Assumption~\ref{ineqforall} holds) and all $(P,t)\in\ptset$. For our running example, we can prove a stronger strict inequality when $(P,t)$ is of the form $(\pab,t)$ :
\[
\nuopt=\max_{k\in\nn} \nu_k<1+\gamma^2(1-\gamma^2)^{-2}
\]
From Example~\ref{exampleprop8}, $H_{\pab,t}(1)=t(a+b)\geq 1+b/a$ and as $b/a> \gamma^2(1-\gamma^2)^{-2}$ we get $H_{\pab,t}(1)> 1+\gamma^2(1-\gamma^2)^{-2}$.

A general statement for this latter stronger strict inequality would be 
$\nuopt<H_{P,t}(1)$ for all $(P,t)$ in the closure of $\ptset$. The set $\ptset$ is not closed in general as $\lyap{A}$ is open in the set of Hermitian matrices.
\qed
\end{example}
We now can construct a formula for the general case similar to the one proposed  in~\cite{DBLP:journals/jota/Adje21} and recalled at Eq.~\eqref{oldformula}. We denote by $\klyap_\nu$ this new formula since it is constructed from solutions of discrete Lyapunov equation. Actually $\klyap_\nu$ depends on the terms of the sequence, the Lyapunov candidate and the scaling to dominate $Q$. Let $k\in\pos_\nu$ and $(P,t)\in \ptset$, we write:
\begin{equation}
\label{Klyapj}
\klyap_\nu(k,P,t):=
\left\lfloor \dfrac{1}{\ln\norm{A}_P}\ln\left(H_{P,t}^{-1}\left(\nu_k\right)\right)\right\rfloor +1\enspace .
\end{equation} 
\begin{theorem}
\label{mainth}
Let $k\in\pos_\nu$ and $(P,t)\in \ptset$. Then :
\begin{enumerate}
\item $\klyap_\nu(k,P,t)$ is well-defined i.e a natural integer;
\item $\klyap_\nu(k,P,t)=\min\{j\in\nn : H_{P,t}(\norm{A}_P^j)< \nu_k\}=\min\{j\in\nn : \norm{A}_P^j< H_{P,t}^{-1}(\nu_k)\}$;
\item $\klyap_\nu(k,P,t)\geq k+1$;
\item For all $j\geq \klyap_\nu(k,P,t)$, $\nu_j<\nu_k$;
\item Let $j\in\pos_\nu$. If $\nu_j\leq \nu_k$ then $\klyap_\nu(k,P,t)\leq \klyap_\nu(j,P,t)$.
\end{enumerate}
\end{theorem}
\begin{proof}
Let $k\in\pos_\nu$ and $(P,t)\in \ptset$.

\noindent ({\itshape 1}) From Prop.~\ref{nukineg} and Assumption~\eqref{ineqforall} and using the fact that $H_{P,t}^{-1}$ is strictly increasing, we have $\ln(H_{P,t}^{-1}(\nu_k))<\ln(H_{P,t}^{-1}(H_{P,t}(1)))=0$. We conclude, using the third point of  Lemma~\ref{lyapnorm} that $\klyap_\nu(k,P,t)\in\nn$.

\noindent ({\itshape 2}) The existence of an integer $j$ s.t $H_{P,t}(\norm{A}_P^j)< \nu_k$ is a consequence of the fact that $\nu_k>0$ and $H_{P,t}(\norm{A}_P^j)$ tends to 0 as $j$ goes to $+\infty$. The equality between the two min follows readily from the fact that $H_{P,t}$ and its inverse are strictly increasing. 

Now, let $j\in\nn$ such that $\norm{A}_P^j< H_{P,t}^{-1}(\nu_k)$. This is equivalent to $j\ln(\norm{A}_P)<\ln(H_{P,t}^{-1}(\nu_k))$ which is the same as $j> \ln(H_{P,t}^{-1}(\nu_k))/\ln(\norm{A}_P)\geq \klyap_\nu(k,P,t)-1$. Since $\ln(\norm{A}_P)<0$, we have $\klyap_\nu(k,P,t)\ln(\norm{A}_P)<\ln(H_{P,t}^{-1}(\nu_k))$ and thus $\norm{A}_P^{\klyap_\nu(k,P,t)}<H_{P,t}^{-1}(\nu_k)$, the result follows.

\noindent ({\itshape 3}) Let $k\in\nn$. From Prop.~\ref{nukineg}, we have $\nu_k\leq H_{P,t}(\norm{A}_P^k)$ which is equivalent, using the natural logarithm, to $\ln(H_{P,t}^{-1}(\nu_k))\leq k\ln(\norm{A}_P)$. As $\ln(\norm{A}_P)<0$ and the integer part function is increasing and stable on the integers, we get $\klyap_\nu(k,P,t)\geq k+1$.

\noindent ({\itshape 4}) Let $j\geq \klyap_\nu(k,P,t)$. As $\norm{A}_P<1$ and by definition of $\klyap_\nu(k,P,t)$, we have $H_{P,t}(\norm{A}_P^j)\leq H_{P,t}(\norm{A}_P^{\klyap_\nu(k,P,t)})<\nu_k$. From Prop.~\ref{nukineg}, $\nu_j\leq H_{P,t}(\norm{A}_P^j)$ which ends the proof.

\noindent ({\itshape 5}) The result follows readily from the fact that the composition of functions $\lfloor\cdot\rfloor \circ (/ \ln(\norm{A}_P)) \circ \ln \circ H_{P,t}^{-1}$ is decreasing on $\rr_+^*$.
\end{proof}
\begin{remark}
The fourth statement of Th~\ref{mainth} is more precise compared with the second assertion of Th~\ref{thold} (from the previous paper). Indeed, in Th~\ref{mainth}, we have stated that the terms of the sequence after $\klyap_\nu(k,P,t)$ are strictly smaller than $\nu_k$ whereas they were only smaller in Th~\ref{thold}. This result permits to gain one iteration in Algorithm~\ref{algoAffQP} and we can replace the stopping condition $k<K$ by $k<K-1$ at Line 17 of Algorithm~\ref{algoAffQP}.
\end{remark}
\begin{proposition}
\label{lyapprop}
If $Q\in\lyap{A}$ and $q=0$. Then $\klyap_\nu(0,Q,1)=1$.
\end{proposition}
\begin{proof}
If $Q\in\lyap{A}$, then $(Q,1)\in\ptset$ and we can define $H_{Q,1}^{-1}$. As $Q\in\lyap{A}$, $Q\succ 0$ and since $q=0$ and $\xin\neq\{0\}$, we have $\nu_0=\mu(Q)>0$ and $0\in\pos_\nu$. Then, we can define $\klyap_\nu(0,Q,1)$. Finally, $\klyap_\nu(0,Q,1)=\left\lfloor \left(\ln\norm{A}_Q\right)^{-1}\ln\left(H_{Q,1}^{-1}\left(\nu_0\right)\right)\right\rfloor +1=\left\lfloor \left(\ln\norm{A}_Q\right)^{-1}\ln\left(\sqrt{\nu_0\left(\mu(Q)\right)^{-1}}\right)\right\rfloor +1=\left\lfloor \left(\ln\norm{A}_Q\right)^{-1}\ln\left(\sqrt{\mu(Q)(\mu(Q))^{-1}}\right)\right\rfloor +1=1$.
\end{proof}
At Prop.~\ref{lyapprop}, we recover the result of Corollary~\ref{corollary1} and the discussion just after. Indeed from the fourth statement of Th~\ref{mainth}, when $Q\in\lyap{A}$ and $q\neq 0$, we have, for all $j\geq 1$, $\nu_j<\nu_0$. This means that the optimal value is only achieved at $k=0$.

At Figure~\ref{drawsimple}, we give a more graphical view of the role played by $\klyap_\nu$. 
Let us consider a reference sequence $({\color{blue}x_k})_k$ to study. Let also consider two sequences, $({\color{red}y_k})_k$ and $({\color{orange} z_k})_k$ s.t.:
\begin{itemize}
\item ${\color{blue}x_k}\leq {\color{red}y_k}$ for all $k\in\nn$;
\item ${\color{blue}x_k}\leq {\color{orange}z_k}$ for all $k\in\nn$;
\item for all $k\in\nn$, ${\color{red}y_k}=h_1(\beta_1^k)$ and ${\color{orange}z_k}=h_2(\beta_2^k)$ where $h_1,h_2$ are homeomorphisms on $\rr_+$ and $\beta_1,\beta_2\in (0,1)$.
\end{itemize}
\begin{figure}[h!]
\begin{center}
\begin{tikzpicture}[scale = 1.4]
 \newcommand{\xint}{0,1,...,56}
    \begin{axis}[
    axis lines=left,
    axis x line =middle,
    ymax=10.7,
    ymin=-2,
    xmax=59,
    minor tick num=1,
    xticklabels = \empty,
    yticklabels = \empty,
    x tick style ={color=black, ultra thick},    
    xtick={5,14,22,28,41,48},
    x label style={at={(1.01,0.07)},anchor=south},
    y label style={at={(-0.01,1.08)},rotate=-90,anchor=north},
    xlabel = {{\scriptsize $k$}},
    ylabel = {{\scriptsize {\color{blue}$x_k$}, {\color{red}$y_k$}, {\color{orange}$z_k$}}},
    axis line style={-{Stealth[scale=1.5]}}
    ]
\foreach \x in \xint {
\addplot[only marks,mark=square,mark size= 0.3mm,blue] coordinates {(\x,{5*exp(-0.0071*(\x-14)^2)})};
}
\label{xk}
\foreach \x in \xint {
\addplot[only marks,mark=diamond,mark size= 0.3mm,red] coordinates {(\x,{16.1*(0.958^(\x))})};
}
\label{yk}
\foreach \x in \xint {
\addplot[only marks,mark=*,mark size= 0.3mm,orange] coordinates {(\x,{8*(0.978^(\x))})};
}
\label{zk}
\draw[thick,dashed,black] (axis cs:5,0) -- (axis cs:5,{5*exp(-0.0071*81});
\draw (axis cs:5,-0.2) node[below,scale=0.4]{$5$};
\draw[thick,dashed,cyan] (axis cs:5,{5*exp(-0.0071*81)}) -- (axis cs:54,{5*exp(-0.0071*81)});
\draw[thick,dashed,gray]  (axis cs:41,{16.1*(0.958^(41))}) -- (axis cs:41,0);
\draw (axis cs:41,-0.2) node[below,scale=0.4]{${\color{red}\klyap_{x,1}(5)}$};
\draw (axis cs:48,-0.2) node[below,scale=0.4]{${\color{orange}\klyap_{x,2}(5)}$};
\draw[thick,dashed,gray]  (axis cs:48,{8*(0.978^(48))}) -- (axis cs:48,0);
\draw[thick,dashed,black] (axis cs:14,0) -- (axis cs:14,{5});
\draw (axis cs:14,-0.2) node[below,scale=0.4]{$14$};
\draw[thick,dashed,cyan] (axis cs:14,{5}) -- (axis cs:34,{5});
\draw[thick,dashed,gray]  (axis cs:28,{16.1*(0.958^(28))}) -- (axis cs:28,0);
\draw (axis cs:21.7,-0.2) node[below,scale=0.4]{${\color{orange}\klyap_{x,2}(14)}$};
\draw (axis cs:28.3,-0.2) node[below,scale=0.4]{${\color{red}\klyap_{x,1}(14)}$};
\draw[thick,dashed,gray]  (axis cs:22,{8*(0.978^(22))}) -- (axis cs:22,0);
\node [draw,fill=white] at (rel axis cs: 0.8,0.8) {\shortstack[l]{
        \ref{xk} {\color{blue}$x_k$} \\
        \ref{yk} {\color{red}$y_k$}\\
        \ref{zk} {\color{orange}$z_k$}}};
\end{axis}
\end{tikzpicture}
\end{center}
\caption{A draw to illustrate Formula~\eqref{Klyapj} and Th~\ref{mainth}}
\label{drawsimple}
\end{figure}
As we analyse the sequence $(x_k)_k$, we are interested in the integers $\klyap_x$. Here, we simplify the notation : $\klyap_{x,1}$ (resp. $\klyap_{x,2}$) means that we use the homeomorphism $h_1$ and thus $(y_k)_k$ (resp. $h_2$ and thus $(z_k)_k$). Moreover, $\klyap_{x,1}$ (resp. $\klyap_{x,2}$) appears in red (resp. in orange). Hence, we have $\klyap_{x,1}(14)=28$, $\klyap_{x,2}(14)=22$,
$\klyap_{x,1}(5)=41$ and $\klyap_{x,2}(5)=48$. The integers computed for $k=14$ are smaller than the one computed for $k=5$ since $x_5< x_{14}$. This illustrates the fifth statement of Th~\ref{mainth}. 

The horizontal dot lines highlight the fourth statement of Th~\ref{mainth}. We see, for example, that for all $j\geq \klyap_{x,1}(5)$, $x_j<x_5$. This integer bound is safe but not tight. Indeed, we see, only using $(x_k)_k$, that $x_j<x_5$ for all $j>24$. Unfortunately, we cannot guess this bound with a close formula. In general, our sequence $(\nu_k)_k$  does not benefit from good mathematical properties (order-preservation, recurrence formula, definition from homeomorphism,...). Thus, we exploit the properties of the upper-bound homeomorphic geometric sequence. With the draw, we see the importance of the computation of  a tight upper-bound. 

Recall that $\bigkse_\nu$ is the smallest rank $k$ such that $\nu_k=\max_{j\in\nn} \nu_j$.
\begin{proposition}
\label{inter}
Let $(P,t)\in \ptset$. The following statements hold:
\begin{enumerate}
\item For all $j\in\pos_\nu$, $\bigkse_\nu<\klyap_\nu(j,P,t)$;
\item We have \[\min_{k\in \pos_\nu} \klyap_\nu(j,P,t)=\klyap_\nu(\bigkse_\nu,P,t)\enspace ;\]
\end{enumerate}
\end{proposition}
\begin{proof}
Let $(P,t)\in \ptset$.

\noindent ({\itshape 1}) Let $j\in\pos_\nu$. If $j<\bigkse_\nu$, then from the fifth point of Th~\ref{mainth}, as $\nu_j<\nuopt=\nu_{\bigkse_\nu}$, we have $\klyap_\nu(\bigkse_\nu,P,t)\leq \klyap_\nu(j,P,t)$. From the first point of Th~\ref{mainth}, $\bigkse_\nu<\bigkse_\nu+1\leq \klyap_\nu(\bigkse_\nu,P,t)\leq \klyap_\nu(j,P,t)$. Now if $j\geq \bigkse_\nu$, from the second statement of Th~\ref{mainth} and by assumption, $\bigkse_\nu\leq j<j+1\leq \klyap_\nu(j,P,t)$
 
 \noindent ({\itshape 2}) Let $k\in\pos_\nu$. If $k\notin\agmx(\nu)$, then $\nu_k<\nuopt=\nu_{\bigkse_\nu}$ and from the fifth statement of Th~\ref{mainth}, $\klyap_\nu(\bigkse_\nu,P,t)\leq \klyap_\nu(k,P,t)$. If $k\in\agmx(\nu)$, $\nu_k=\nuopt=\nu_{\bigkse_\nu}$ and $\klyap_\nu(\bigkse_\nu,P,t)=\klyap_\nu(k,P,t)$. 
 \end{proof}
\begin{example}
\label{krunning}
Applying the computation to our running example leads, for all $k\in\nn$, and for all $a,b\in\lyapex(\gamma)$ and all $t\in\cont(\pab)$, to:
\[
\begin{array}{l}
\klyap_\nu(k,\pab,t)-1\\
=\left\lfloor \dfrac{1}{\ln\norm{A}_{\pab}}\ln\left(\dfrac{\sqrt{4t\nu_k}}{2t\sqrt{\mu(\pab)}}\right)\right\rfloor
=\left\lfloor\dfrac{\ln\left(\dfrac{\gamma^{2k}(1+k)^2}{t(a+b)}\right)}{\ln\left(\dfrac{\gamma^2}{2}\left(\dfrac{a+2b+\sqrt{a(a+4b)}}{b}\right)\right)}\right\rfloor\\
=\left\lfloor\dfrac{\ln\left(\dfrac{\gamma^{2k}(1+k)^2}{1+\frac{b}{a}}\right)-\ln(at)}{\ln\left(\dfrac{\gamma^2}{2}\left(\dfrac{1+2\frac{b}{a}+\sqrt{1+4\frac{b}{a}}}{\frac{b}{a}}\right)\right)}\right\rfloor
\end{array}
\]
Applying the same change of variable as earlier $b'=b/a$, we obtain:
\[
\begin{array}{l}
\left\lfloor\dfrac{\ln\left(\dfrac{\gamma^{2k}(1+k)^2}{1+\frac{b}{a}}\right)-\ln(ta)}{\ln\left(\dfrac{\gamma^2}{2}\left(\dfrac{1+2\frac{b}{a}+\sqrt{1+4\frac{b}{a}}}{\frac{b}{a}}\right)\right)}\right\rfloor
=\left\lfloor\dfrac{\ln\left(\dfrac{\gamma^{2k}(1+k)^2}{1+b'}\right)-\ln(ta)}{\ln\left(\dfrac{\gamma^2}{2}\left(\dfrac{1+2b'+\sqrt{1+4b'}}{b'}\right)\right)}\right\rfloor
\end{array}
\]
The latter function increases as $ta$ increases since the denominator is negative. We have already shown at Example~\ref{exampleprop8} that $ta\geq 1$ for all $a,b\in\rr$ and $t$ s.t. $(\pab,t)\in\ptex$. The best choice is to have $\ln(at)=0$. Finally again we reduce the choice of an element $\pab\in\lyap{A_\gamma}$ to a matrix of the form $\pub$ :
\[
\klyap_\nu(k,P_{1,b},1)-1=\left\lfloor\left(\ln\left(\dfrac{\gamma^2}{2}\left(\dfrac{1+2b+\sqrt{1+4b)}}{b}\right)\right)\right)^{-1}\ln\left(\dfrac{\gamma^{2k}(1+k)^2}{(1+b)}\right)\right\rfloor
\]
Now we consider numerical values and let us take, for example, $\gamma=\exp(-1/20)$. In this case, as mentioned at Subsection~\ref{runningex}, the integer $k_1^*=19$ is optimal and provides an optimal value equals to $400\exp(-1.9)\simeq 59.827$. Following Lyapunov conditions i.e. Eq~\eqref{lyaprunn}, we have to choose $b>\exp(-1/10)/((1-\exp(-1/10))^2)\simeq 99.917$. In consequence, let us take first $b=100$. We complete the table by the values of $\klyap_\nu(k,P_{1,b},1)$ for $b=500$ and $b=1000$ to analyze the behaviour of $\klyap_\nu(k,P_{1,b},1)$ with respect to $b$.
\begin{center}
\begin{table}[h!]
{\footnotesize
\begin{center}
\begin{tabular}{|*{9}{c|}}
\hline
k& 0 & 1& 2& 3& 10 & {\bf 19} & 25 & 100 \\
\hline
$\nu_k$ & 1 & 3.6193 &7.3686& 11.853 & 44.513 & {\bf 59.827} & 55.489  & 0.4631\\
\hline 
$b=100$& 110888 & 79882 & 62901 & 51479 &19687&{\bf 12582} & 14391 & 129383\\
\hline
$b=500$ & 113 & 90 & 77 & 68 &44&{\bf 39} & 40 & 127\\
\hline 
$b=1000$& 102 & 83 & 72 & 65 &46&{\bf 42} & 43 & 113\\
\hline
\end{tabular}
\end{center}
}
\caption{Integers $\klyap_\nu(k,P_{1,b},1)$ for $\gamma=\exp(-1/20)$ and $b=100,500$ and $1000$}
\label{firsttable}
\end{table}
\end{center}
We recall that the integers $\klyap_\nu(k,P_{1,b},1)$ are used to limit the number of quadratic maximization problems to solve to compute $\nuopt$. In Algorithm~\ref{algoAffQP}, as we cannot guess $\kopt$, we compute the sequence of $(\klyap_\nu(k,P_{1,b},1)-1)_{k\in\Gamma(\nu)}$ as maximal iteration numbers for the loop at Line 17. This maximal iteration number is not fixed and decreases at each $k\in\Gamma(\nu)$. The best stopping condition is achieved at $k=\bigkse_\nu$. However, we cannot identify an integer $k$ as the optimal solution $\kopt=\bigkse_\nu$ without computing the rest of the sequence until the rank $\klyap_\nu(\bigkse_\nu,P_{1,b},1)-1$. Hence, for this numerical example, the integers $\klyap_\nu(k,P_{1,b},1)$ after $k=19$ are neither computed in practice nor used in Algorithm~\ref{algoAffQP}. Furthermore, we compute $\nu_k$ until $k=\klyap_\nu(19)-1=38$. 

At Table~\ref{firsttable}, we observe the third and the fifth statements of Th~\ref{mainth} in the evolution of $ \klyap_\nu(k,P_{1,b},1)$ with respect to $k$. The integers $\klyap_\nu(k,P_{1,b},1)$ are greater than $k+1$ and when $\nu_j\leq \nu_k$, $\klyap_\nu(k,P_{1,b},1)\leq \klyap_\nu(j,P_{1,b},1)$. A good choice of $b$ (or more generally of $P\in\lyap{A}$) is important to reduce the number of quadratic maximization problems to solve. The integers of the line $b=100$ are very problematic. This is caused by the fact that $\norm{A_\gamma}_{P_{1,b}}^2$ is very close to 1 ($\simeq 0.999959$). For $b=500$ and $b=1000$, $\norm{A_\gamma}_{P_{1,b}}^2$ are respectively close to $0.946218$ and $0.933907$. The square of the spectral radius of $A_\gamma$ is equal to $\exp(-0.1)\simeq 0.904837$. Finally, at Table~\ref{firsttable}, we see that $b$ may not be optimal everywhere. For example, $b=1000$ is better at $k=0,1,2,3$ and $100$ but worse at $k=19$. 
\qed  
\end{example}
In Prop.~\ref{inter}, we minimized $\klyap_\nu(k,P,t)$ over integers $k$ and found that it was achieved at $k=\bigkse_\nu$ i.e. the first rank that achieves $\max_{j\in\nn} \nu_j$. Now, we propose to minimize $\klyap_\nu(k,P,t)$ with respect to $(P,t)\in\ptset$ when $k$ is fixed. In Algorithm~\ref{algoAffQP}, it allows to change the Lyapunov candidate and the scaling factor $t$ at every $k$ in $\Gamma(\nu)$. The only reduction that can be done is on $t$ using the fact that $t\mapsto \klyap_\nu(k,P,t)$ is increasing. This fact is at Example~\ref{krunning} to rewrite $\klyap_\nu(k,\pub,1)$.
\begin{theorem}
\label{thfond}
Let $k\in\pos_\nu$. Let us define :
\[
\bigklyap_\nu(k):=\displaystyle{\min_{(P,t)\in\ptset} \klyap_\nu(k,P,t)}
\enspace .
\] 
Then :
\begin{enumerate}
\item For all $P\in\lyap{A}$, the function $\cont(P)\ni t\mapsto \klyap_\nu(j,P,t)$ is increasing;
\item We have:
\[
\bigklyap_\nu(k)=\min_{P\in\lyap{A}} \klyap_\nu\left(k,P,\topt{P}\right)
\] 
\item $\bigkse_\nu< \bigklyap_\nu(k)$.
\item $\displaystyle{\bigklyap_\nu(\bigkse_\nu)=\min_{j\in\pos_\nu} \bigklyap_\nu(j)}$;
\item For all $k\in\pos_\nu$, 
\[    
    \max_{j\in\nn} \nu_j =\max_{0\leq j \leq \bigklyap_\nu(\bigkse_\nu)-1} \nu_j=\max_{0\leq j\leq \klyap_\nu(k)-1} \nu_j
    \]
\end{enumerate}
\end{theorem}
We recall that for all $(P,t)\in\ptset$, $\mu(P)>0$ and $t$ cannot be null when $Q\npreceq 0$. On the other hand, when $Q\preceq 0$, $t$ can be null but we have assumed that $q\neq 0$. Finally, the real $(\sqrt{4tx+(\norm{q}_P^*)^2}+\norm{q}_P^*)\sqrt{\mu(P)}$ is not null for all $x>0$.
\begin{lemma}
\label{lemmaprmin}
Let us introduce, for all $P\succ 0$ and $t\in\rr_+$, the function on $\rr_+^*$:
\[
x\mapsto h_{P,t}^{-}(x):=\dfrac{2x}{(\sqrt{4tx+(\norm{q}_P^*)^2}+\norm{q}_P^*)\sqrt{\mu(P)}}
\]
Then:
\begin{enumerate}
\item for all $(P,t)\in\ptset$, $H_{P,t}^{-1}=h_{P,t}^{-}$;
\item for all $x>0$, $u\mapsto h_{P,u}^{-}(x)$ is decreasing on $\rr_+$.
\end{enumerate}
\end{lemma}
%
\begin{proof}[Proof of Theorem~\ref{thfond}]
\noindent ({\itshape 1}) The result follows from Lemma~\ref{lemmaprmin} and the fact that $\norm{A}<1$. 

\noindent ({\itshape 2}) This is a direct consequence of the first point and Lemma~\ref{propcont}.

\noindent ({\itshape 3}) From Prop.~\ref{inter}, for all $(P,t)\in\ptset$, for all $k\in\pos_\nu$, $\bigkse_\nu+1\leq 
\klyap_\nu(k,P,t)$ and then $\bigkse_\nu+1\leq 
\displaystyle{\min_{(P,t)\in\ptset} \klyap_\nu(k,P,t)}=\klyap_\nu(k)$. 

\noindent ({\itshape 4}) As independent infima commute: 
\[
\begin{array}{ll}
\displaystyle{\min_{k\in\pos_\nu} \bigklyap_\nu(k)}
&=\displaystyle{\min_{k\in\pos_\nu}\min_{(P,t)\in\ptset} \klyap_\nu(k,P,t)}\\
&=\displaystyle{\min_{(P,t)\in\ptset} \min_{k\in\pos_\nu}\klyap_\nu(k,P,t)}
\end{array}
\enspace .
\]
From Prop.~\ref{inter} we conclude that $\displaystyle{\min_{k\in\pos_\nu} \bigklyap_\nu(k)}=\min_{(P,t)\in\ptset} \klyap_\nu(\bigkse_\nu,P,t)=\bigklyap_\nu(\bigkse_\nu)$.

\noindent({\itshape 5}) The results follows readily from the third statement of Prop.~\ref{supremumpos} as $\bigklyap_\nu(\bigkse_\nu)$ and 
for all $k\in\pos_\nu$, $\klyap_\nu(k)$ are strictly greater than $\bigkse_\nu$. 
\end{proof}
Let us define for all $k\in\pos_\nu$, for all $P\in\lyap{A}$:
\begin{equation}
\label{realobj}
G_k(P):=\dfrac{\ln\left(H_{P,t_Q^*(P)}(\nu_k)\right)}{\ln\left(\norm{A}_P\right)}\quad \text{ and }\quad  \lyopt_k=\inf\left\{G_k(P) : P\in\lyap{A}\right\}
\end{equation}
We have, for all $k\in\pos_\nu$, for all $P\in\lyap{A}$, $\klyap_\nu(k,P,\topt{P})=\lfloor G_k(P)\rfloor+1$. Moreover,$\lyopt_k$ is finite since $G_k$ is greater than $\bigkse_\nu$ on $\lyap{A}$.
\begin{proposition}
\label{propopt}
Let $k\in\pos_\nu$. There exists $\overline{P_k}\in\lyap{A}$ such that $\bigklyap_\nu(k)=\lfloor G_k(\overline{P_k})\rfloor+1$ and $\bigklyap_\nu(k)=\lfloor \lyopt_k\rfloor+1$.
\end{proposition}
\begin{proof}
Let $k\in\pos_\nu$. Since $\{\klyap_\nu(k,P,\topt{P}) : P\in\lyap{A}\}$ is lower-bounded set of integers, it has a smallest element achieved at some  $\overline{P_k}\in\lyap{A}$. As $\bigklyap_\nu(k)-1=\lfloor G_k(\overline{P_k})\rfloor$ and $\lyopt_k \leq G_k(\overline{P_k})$, we conclude that $\lfloor\lyopt_k\rfloor+1 \leq \bigklyap_\nu(k)$. Conversely, for all $P\in\lyap{A}$, $\lfloor G_k(P)\rfloor\leq G_k(P)$ which implies that
$\bigklyap_\nu(k)\leq \lyopt_k+1$ and thus $\bigklyap_\nu(k)\leq \lfloor \lyopt_k\rfloor +1$.
\end{proof}
We apply Theorem~\ref{thfond} and Prop.~\ref{propopt} to our running example to reduce the values $\klyap_\nu(k,P,t)$ computed at Table~\ref{firsttable}. Note that in our running example, we have restricted the Lyapunov candidates to diagonal elements of $\lyap{A_\gamma}$. In this case, as we shown at Example~\ref{krunning}, that it suffices to minimize $\klyap_\nu(k,\pub,1)$ over $b>\gamma^2(1-\gamma^2)^{-2}$. This is a one-dimensional minimization problem for which we can apply a dichotomy method.
\begin{example}
\label{minexamp}
Recall that for $k\in\nn$, $\nu_k=\gamma^{2k}(1+k)^2$. In our running example, to minimize $b\mapsto \klyap_\nu(k,\pub,1)$ on $(\gamma^2(1-\gamma^2)^{-2},+\infty)$, we minimize $b\mapsto G_k(\pub)$ on $(\gamma^2(1-\gamma^2)^{-2},+\infty)$ which is a nonlinear program on the variable $b$:
\[
\Min_b\ \left(\ln\left(\dfrac{\gamma^2}{2}\left(\dfrac{1+2b+\sqrt{1+4b)}}{b}\right)\right)\right)^{-1}\ln\left(\dfrac{\nu_k}{(1+b)}\right)\quad \st \quad
b>\dfrac{\gamma^2}{(1-\gamma^2)^2} 
\]
The objective function is, in general, not convex as we can see at Figure~\ref{nonconvex} for $\gamma=0.9$ and $k=10$. For example, in this case, we have $2G_k(P_{1,500})\simeq 42.505>G_k(P_{1,200})+G_k(P_{1,800})\simeq 41.467$.

In Example~\ref{complicated}, we wrote $B=(1+2b+\sqrt{1+4b})(2b)^{-1}$ which was proved to be in $(1,\gamma^{-2})$. This leads to the nonlinear program:
\[
\Min_B\ \dfrac{1}{\ln(B\gamma^2)}\ln\left(\dfrac{(B-1)^2\nu_k}{B^2-B+1}\right)\quad \st \quad 1<B<\gamma^{-2}
\]
The second formulation brings finite bounds around the decision variable $B$. Then we can solve the nonlinear program using the dichotomy method applied on the derivative of \[g_k:B\mapsto \dfrac{1}{\ln(B\gamma^2)}\ln\left(\dfrac{(B-1)^2\nu_k}{B^2-B+1}\right)\enspace.\] Indeed:
\begin{enumerate}
\item The nonlinear equation $g'_k(B)=0$ has a solution on $(1,\gamma^{-2})$;
\item For all $\overline{B}\in(1,\gamma^{-2})$ such that $g'_k(\overline{B})=0$, we have ${g_k''}(\overline{B})>0$. In other words, all critical points are local minimizers. We conclude that there exists one critical point which is the unique minimizer of $g_k$;
\item The dichotomy method converges to this unique minimizer.
\end{enumerate}
A proof for this result is provided in Appendix.
\begin{figure}[h!]
\begin{center}
\begin{tikzpicture}[scale = 1.0]
 \begin{axis}[
    axis lines=left,
    xmin=0,
    xmax=1250,
    ymin=17,
    ymax=24.2,
    x label style={at={(1.01,0)},anchor=south},
    y label style={at={(0,1.1)},rotate=-90,anchor=north},
    xlabel = {$b$},
    ylabel = {$G_k(P_{1,b})$},
   samples=2000,domain=25:2000,restrict y to domain =17:24]
    \addplot[very thick,orange] plot (\x, {ln(121*0.9^(20)/(1+\x))/ln(0.5*0.81*(1+2*\x+sqrt(1+4*\x))/\x))});
\end{axis}
\end{tikzpicture}
\end{center}
\caption{The graph of $G_k(P_{1,b})$ for $\gamma=0.9$ and $k=10$}
\label{nonconvex}
\end{figure}

Once the minimizer $\overline{B}$ computed, we compute $\overline{b}=\oB (\oB-1)^{-2}$ and $\klyap_\nu(k,P_{1,\overline{b}},1)=\lfloor G_k(P_{1,\overline{b}})\rfloor+1$.

We come back to the case where $\gamma=\exp(-1/20)$. When we apply the dichotomy method to $g_k$ we obtain the integers $\klyap_\nu(k,P_{1,b},1)$ of Table~\ref{finaltable}. We compare two series of integers $\klyap_\nu(k,P_{1,b},1)$. The integers of the line $\klyap_\nu(k,P_{1,b_k},1)$ are computed from the minimizer  $B_k$ of $g_k$ at each $k$ whereas the integers $\klyap_\nu(k,P_{1,b_0},1)$ are computed from  $B_0$, the minimizer of $g_0$. This comparison aims to analyze the gain obtained by an optimal change of $\pub$ at each $k$.   
\begin{center}
\begin{table}[h!]
{\footnotesize
\begin{center}
\begin{tabular}{|*{8}{c|}}
\hline
$k$& 0 & 1& 3& 10 & {\bf 19} & 25 & 100 \\
\hline
$\nu_k$ & 1 & 3.6193 & 11.853 & 44.513 & {\bf 59.827} & 55.489  & 0.4631\\
\hline
$B_k$& 
 1.02065& 
   1.0249 & 
   1.0313& 1.0466 & {\bf 1.0539}  & 1.0517 & 1.0188 \\
\hline
$b_k$&  2392.9& 1651.1  & 1053.7 & 481 & {\bf 362.7} & 392.8 & 2883.4\\
\hline
$\klyap_\nu(k,P_{1,b_k},1)$&  98  & 82  & 65& 44 & {\bf 38} & 40  & 108\\
\hline
$\klyap_\nu(k,P_{1,b_0},1)$ &  98  & 83 & 67& 51 & {\bf 47} & 48  & 108\\
\hline
\end{tabular}
\end{center}
}
\caption{The integers $\klyap_\nu(k,P_{1,b},1)$ for an optimized $b$ or not and $\gamma=\exp(-1/20)$}
\label{finaltable}
\end{table}
\end{center}
The differences between the rows $\klyap_\nu(k,P_{1,b_k},1)$ and $\klyap_\nu(k,P_{1,b_0},1)$ seem to be unsignificant except near the optimal $k=\bigkse_\nu$. The dynamic update of $\pub$ reduces the number of iterations by 9 and the price to pay is the use of a dichotomy method at each $k$. 
\qed
\end{example}
\subsection{Running Example with Numerical Values}
We propose to do the same comparisons for the running example but for different values of $\gamma\in(0,1)$. The implementations have been done using Octave 6.3.0~\cite{octave}.

\paragraph*{The case where $\gamma=1/3$.} Since $\gamma<\exp(-1)$, we have $\kopt=0$ and $\nuopt=1$. The different integers $\klyap_\nu(k,\pub,1)$ are given at Table~\ref{tableuntiers}.
\begin{center}
\begin{table}[h!]
{\footnotesize
\begin{center}
\begin{tabular}{|*{7}{c|}}
\hline
\multicolumn{7}{|c|}{$\gamma=1/3$}\\
\hline
$k$& {\bf 0} & 1& 2& 3& 10  & 100 \\
\hline
$\nu_k$ & {\bf 1} & 0.4444 &0.111& 2.19e-02 & 3.47e-08  & 3.84e-92\\
\hline
$\klyap_\nu(k,P_{1,b_k},1)$& {\bf 1} & 2  & 3  & 4& 11 & 101\\
\hline
$\klyap_\nu(k,P_{1,b_0},1)$&  {\bf 1}  & 2 & 4  & 6& 24 & 288\\
\hline
\end{tabular}
\end{center}
}
\caption{The integers $\klyap_\nu(k,P_{1,b},1)$ for an optimized $b$ or not $\gamma=1/3$}
\label{tableuntiers}
\end{table}
\end{center}
The integers computed after $k=0$ are not used in practice since $\klyap_\nu(0,P_{1,b_0},1)=1$.
\paragraph*{The case where $\gamma=0.5$.}
As $\gamma>\exp(-1)$ and $\lfloor -(\ln(\gamma)^{-1}-1\rfloor=(1-\gamma)^{-1}-2=0$, $k=0$ and $k=1$ are both optimal and $\nuopt=1$. The different integers $\klyap_\nu(k,\pub,1)$ are given at Table~\ref{tableundemi}.
\begin{center}
\begin{table}[h!]
{\footnotesize
\begin{center}
\begin{tabular}{|*{7}{c|}}
\hline
\multicolumn{7}{|c|}{$\gamma=0.5$}\\
\hline
$k$& {\bf 0} & 1& 2& 3& 10  & 100 \\
\hline
$\nu_k$ & {\bf 1} & 1 &0.5625& 0.25 & 1.15e-04  & 6.34e-57\\
\hline
$\klyap_\nu(k,P_{1,b_k},1)$& {\bf 2} & 2  & 3  & 4& 11 & 101\\
\hline
$\klyap_\nu(k,P_{1,b_0},1)$ &  {\bf 2}  & 2 & 3  & 4& 18 & 225\\
\hline
\end{tabular}
\end{center}
}
\caption{The integers $\klyap_\nu(k,P_{1,b},1)$ for an optimized $b$ or not and $\gamma=1/2$}
\label{tableundemi}
\end{table}
\end{center}
Again, the integers $\klyap_\nu(k,P_{1,b},1)$ after $k>0$ are given for the example besides their uselessness. In practice, $\klyap_\nu(0,P_{1,b},1)=2$ implies that we only compute $\nu_0$ and $\nu_1$ to obtain $\nuopt$.

\paragraph*{The family $\gamma_n=(n+0.1)/(n+1.1)$.}
We study the family for $n=20$ or $n=100$ i.e. $\gamma_{20}=201/211$ or $\gamma_{100}=1001/1011$. We recall that, for this family, the only maximizer is  $k_2^*=n$. Hence, for $n=20$, the maximizer is $20$ whereas for $n=100$, it is $100$. The different integers $\klyap_\nu(k,\pub,1)$ are given at Table~\ref{tablensurn}.
\begin{center}
\begin{table}[h!]
{\footnotesize
\begin{center}
\begin{tabular}{|*{8}{c|}}
\hline
\multicolumn{8}{|c|}{$\gamma=201/211$}\\
\hline
$k$&  0 & 1& 2 & 3 & 5 & 10 & {\bf 20} \\
\hline
$\nu_k$ & 1 & 3.6299  & 7.4113 & 11.956 & 22.153  &45.820 & {\bf 63.239}\\
\hline
$\klyap_\nu(k,P_{1,b_k},1)$&  102 & 85 & 75 & 68 & 59  &46 & {\bf 40}\\
\hline
$\klyap_\nu(k,P_{1,b_0},1)$ &  102  & 85 & 76 & 70 & 62 & 53 & {\bf 48}\\
\hline
\hline
\multicolumn{8}{|c|}{$\gamma=1001/1011$}\\
\hline
$k$&  0 & 1& 5 & 20 & 50  & 80 & {\bf 100} \\
\hline
$\nu_k$ & 1 & 3.9213 & 32.594 & 296.32 & 962.57 & 1337.3 & {\bf 1397.1}\\
\hline
$\klyap_\nu(k,P_{1,b_k},1)$&  689  & 608 & 476 & 328 & 235 & 204 & {\bf 200}\\
\hline
$\klyap_\nu(k,P_{1,b_0},1)$&  689  &608 & 483 & 353 & 284 & 265 & {\bf 262}\\
\hline
\end{tabular}
\end{center}
}
\caption{The integers $\klyap_\nu(k,P_{1,b},1)$ for an optimized $b$ or not for $\gamma_{20}=201/211$ and $\gamma_{100}=1001/1011$}
\label{tablensurn}
\end{table}
\end{center}

\paragraph*{The family $\gamma_n=1-10^{-n}$.}
We study the family for $n=2$ or $n=3$ i.e. $\gamma_2=0.99$ or $\gamma_3=0.999$. We recall that, for this family, there exist two maximizers $10^n-2$ and $10^n-1$. Hence, for $n=2$, the maximizers are $98$ and $99$ whereas for $n=3$, they are $998$ and $999$. The different integers $\klyap_\nu(k,\pub,1)$ are given at Table~\ref{table999}.
\begin{center}
\begin{table}[h!]
{\footnotesize
\begin{center}
\begin{tabular}{|*{8}{c|}}
\hline
\multicolumn{8}{|c|}{$\gamma_2=0.99$}\\
\hline
$k$&  0 & 1& 10 & 20 & 50 & 80 & {\bf 98} \\
\hline
$\nu_k$ & 1 & 3.9204  & 98.967 & 295.02 & 952.05  &1314.0 & {\bf 1367.0}\\
\hline
$\klyap_\nu(k,P_{1,b_k},1)$&  680 & 600   & 398 & 323 & 231  &201 & {\bf 198}\\
\hline
$\klyap_\nu(k,P_{1,b_0},1)$ &  680  & 600  & 412 & 348 & 280 & 261 & {\bf 259}\\
\hline
\hline
\multicolumn{8}{|c|}{$\gamma_3=0.999$}\\
\hline
$k$&  0 & 10& 100 & 200 & 500 & 800 & {\bf 998} \\
\hline
$\nu_k$ & 1 & 118.60 & 8351.0 & 27076 & 92292  &129433 & {\bf 135470}\\
\hline
$\klyap_\nu(k,P_{1,b_k},1)$&  9457  & 6731 & 4112 & 3306 & 2353 & 2043 & {\bf 1998}\\
\hline
$\klyap_\nu(k,P_{1,b_0},1)$&  9457  & 6788& 4411 & 3754 & 3068 & 2879 & {\bf 2854}\\
\hline
\end{tabular}
\end{center}
}
\caption{The integers $\klyap_\nu(k,P_{1,b},1)$ for an optimized $b$ or not for $\gamma_2=0.99$ and $\gamma_3=0.999$}
\label{table999}
\end{table}
\end{center}

\paragraph*{The family $\gamma_n=1-10^{-n}+10^{-n-1}$}
Again, we study the family for two values of $n$; for $n=2$ or $n=3$. This means that $\gamma_2=0.991$ or $\gamma_3=0.9991$.
We recall that for this family, the only maximizer is $\sum_{i=1}^n 10^i$. Hence, for $n=2$, the maximizer is $110$ whereas for $n=3$, it is $1110$. The different integers $\klyap_\nu(k,\pub,1)$ are given at Table~\ref{table9991}.
\begin{center}
\begin{table}[h!]
{\footnotesize
\begin{center}
\begin{tabular}{|*{7}{c|}}
\hline
\multicolumn{7}{|c|}{$\gamma=0.991$}\\
\hline
$k$&  0 & 1& 10 & 20 & 80 & {\bf 110} \\
\hline
$\nu_k$ & 1 & 3.9283  & 100.99& 307.17  &1544.4 & {\bf 1686.0}\\
\hline
$\klyap_\nu(k,P_{1,b_k},1)$&  769 & 681  & 456 & 372 &229 & {\bf 220}\\
\hline
$\klyap_\nu(k,P_{1,b_0},1)$ &  769  & 681 & 471 & 391 &  295 & {\bf 289}\\
\hline
\hline
\multicolumn{7}{|c|}{$\gamma=0.9991$}\\
\hline
$k$&  0 & 10& 100 & 500 & 1000 & {\bf 1110} \\
\hline
$\nu_k$ & 1 & 118.84&8519.9& 102008.04  &165495.46 & {\bf 167231.09} \\
\hline
$\klyap_\nu(k,P_{1,b_k},1)$&  10639 & 7615 & 4708  & 2720 & 2231  & {\bf 2220} \\
\hline
$\klyap_\nu(k,P_{1,b_0},1)$& 10639 & 7677 & 5028  & 3488 & 3188 & {\bf 3182} \\
\hline
\end{tabular}
\end{center}
}
\caption{The integers $\klyap_\nu(k,P_{1,b},1)$ for an optimized $b$ or not for $\gamma_2=0.991$ and $\gamma_3=0.9991$}
\label{table9991}
\end{table}
\end{center}
\section{On the Computation of a Solution of the Lyapunov Equation for the Key Formula}
\label{minlyapunov}
\subsection{Nonlinear Semidefinite Programs}
\label{nonlinearSDP}
In this subsection, we study the minimization problem  introduced at Eq~\eqref{realobj}. For the running example, we have restricted $\lyap{A}$ to the subset composed of diagonal elements of it. This restriction allows us to use a dichotomy method to minimize  $b\mapsto G_k(\pub)$ on $(\gamma^2(1-\gamma^2)^{-2},+\infty)$. In general, those simplifications cannot be done. The objective of this subsection is to explain how to rewrite $\inf\{ G_k(P): P\in\lyap{A}\}$ as an optimization problem for which an optimal solution is proved to exist. The optimization problem depicted at ~\eqref{realobj} has the disadvantage to have an open constraints set. Classically, because theoretical and numerical aspects, continuous optimization problems have closed constraints set. However, the closure of $\lyap{A}$ contains matrices $P$ such that $\norm{A}_P=1$ and for which $G_k$ is not defined. Then, we have to consider compact subsets of $\lyap{A}$ for which the optimal solution is proved to exist and is also an optimal solution for $\min\{ \klyap_\nu(k,P,\topt{P}) : P\in\lyap{A}\}$. 

In Prop.~\ref{propoptim}, we propose to bound the constraint set using an homogeneous function $f$ (for example $\lambda_{\rm max}$).  
However, the constraint sets of the minimization problems of Prop.~\ref{propoptim} are not closed and may not have optimal solutions.
\begin{proposition}
\label{propoptim}
Let $f:\hddp\mapsto \rr_+^*$ be an homogeneous function and $\beta$ a strictly positive scalar. Let us define :
\[
\begin{array}{ll}
& V_k^{f,\beta}:=\inf\left\{G_k(P) :P\in\lyap{A}\cap f^{-1}\left((0,\beta]\right)\right\}\\
\text{ and } & W_k^{f,\beta}:=\inf\left\{G_k(P) :P\in\lyap{A}\cap f^{-1}\left(\{\beta\}\right)\right\}
\end{array}
\]
Then $\lyopt_k=V_k^{f,\beta}=W_k^{f,\beta}$.
\end{proposition}
\begin{proof}
Let $P\in\lyap{A}$. Recall that by homogeneity of $\lmax{\cdot}$, we have for all $\alpha>0$, $t_Q^*(\alpha P)=\alpha^{-1} t_Q^*(P)$. Let us define $P'=\beta P/f(P)$. Now, since $(P,t_Q^*(P))\in\ptset$ (and thus also $(P',t_Q^*(P'))$) and $\beta$ and $f(P)$ are strictly positive, we have from the fourth statement of Prop.~\ref{homeoprop}, $H_{P',t_Q^*(P')}^{-1}(\nu_k)=H_{P,t_Q^*(P)}^{-1}(\nu_k)$. Finally, since $P\mapsto \norm{A}_P$ is invariant under scalings, we conclude that $\norm{A}_{P'}=\norm{A}_P$ and then $G_k(P')=G_k(P)$.

We have obviously $\lyopt_k\leq V_k^{f,\beta}\leq W_k^{f,\beta}$. Now, suppose that $\lyopt_k<W_k^{f,\beta}$. Then there exists $R\in\lyap{A}$ such that $G_k(R)<W_k^{f,\beta}$. But $R'=\beta R/f(R)$ belongs to $\lyap{A}\cap f^{-1}\left(\{\beta\}\right)$ since $\lyap{A}$ is a cone and $f$ is homogeneous. Moreover, $W_k^f\leq G_k(R')=G_k(R)$ which contradicts the fact that $G_k(R)<W_k^{f,\beta}$.
\end{proof}
Let us define for all $k\in\pos_\nu$, for all $P\in\lyap{A}$:
\[
\varepsilon_k(P):=\left\lfloor G_k(P)\right\rfloor+1-G_k(P)
\]
\begin{proposition}
\label{optimerror}
Let $\overline{P}\in\lyap{A}$, $\varepsilon<\varepsilon_k(\overline{P})$ and $P\in\lyap{A}$. If  $G_k(\overline{P})\leq G_k(P)\leq G_k(\overline{P})+\varepsilon$ then $\lfloor G_k(P)\rfloor=\lfloor G_k(\overline{P})\rfloor$.
\end{proposition}
\begin{proof}
The result follows from the fact that the integer part function is increasing and $G_k(P)$ is strictly smaller than $\lfloor G_k(\overline{P})\rfloor+1$.
\end{proof}
Prop~\ref{optimerror} asserts that, as we are interested in integer part of $G_k(P)$, we can allow small approximations to compute this integer part. This observation permits to consider $\varepsilon$-optimal solutions in minimization problems of Prop.~\ref{propoptim}. Now, we approximate the constraints set of Prop.~\ref{propoptim} by compact sets representable in machine. 
\begin{lemma}
\label{mori}
Let $\alpha>0$. The unique positive definite matrix $P$ such that $P-A^\intercal P A=\alpha\Idd$ satisfies : 
\[
P\preceq \alpha\lmax{\sum_{j=0}^{d-1} \left(A^\intercal\right)^j A^j} R_I
\]
where $R_I$ is the unique positive definite matrix satisfying $R_I-(A^\intercal)^d R_I A^d=\Idd$.   
\end{lemma}
Lemma~\ref{mori} is a straightforward application of~\cite[Th. 1]{mori1985} for which $B=\sqrt{\alpha}\Idd$ and thus $M=\sqrt{\alpha} \left(\Idd, A^\intercal, (A^\intercal)^2,\ldots,(A^\intercal)^{d-1}\right)$ .

Let us define the following subset of $\hddp$:
\begin{equation}
\label{setsdp}
\sdp{\alpha,\beta}:=\{P\in\hp : P-A^\intercal P A\succeq \alpha\Idd \text{ and }P\preceq \beta \Idd\}
\end{equation}
and the functions on $\rr_+$:
\begin{equation}
\begin{array}{ll}
&\alpha \mapsto \betaf(\alpha):=\alpha \lmax{\sum_{j=0}^{d-1} \left(A^\intercal\right)^j A^j} \lmax{R_I}\\
\\
\text{and} & \beta \mapsto \alphaf(\beta):=\dfrac{\beta}{ \lmax{\sum_{j=0}^{d-1} \left(A^\intercal\right)^j A^j} \lmax{R_I})}
\end{array}
\end{equation}
We observe that $\betaf(\alpha)$ is a value which allows to upper bound $P$ of Lemma~\ref{mori} by $\betaf(\alpha)\Idd$ whereas the value $\alphaf(\beta)$ ensures that $P-A^\intercal P A=\alphaf(\beta)\Idd$ and $P\preceq \beta \Idd$. 
\begin{proposition}
\label{propmori}
The values $\lmax{\sum_{j=0}^{d-1} \left(A^\intercal\right)^j A^j}$ and $\lmax{R_I}$ are strictly greater than 1. Hence $\alphaf$ is well-defined, for all $\beta>0$, $\alphaf(\beta)<\beta$ and for all $\alpha>0$, $\betaf(\alpha)\geq \alpha$. Moreover, for all $\alpha,\beta>0$, $\alpha\leq\alphaf(\beta)\iff \betaf(\alpha)\leq \beta$.  
\end{proposition}
\begin{proof}
The Weyl's inequalities imply that $\lmax{R_I}\geq \lmin{R_I-(A^\intercal)^d R_I A^d}+\lmax{(A^\intercal)^d R_I A^d}\geq 1+\lmax{(A^\intercal)^d R_I A^d}$. Since $R_I$ is positive definite, the matrix $(A^d)^\intercal R_I A^d$ is positive semidefinite. Hence, $\lmax{(A^d)^\intercal R_I A^d}=0$ if and only if $x^\intercal (A^d) ^\intercal R_I A^d x=0$ for all $x\in\rd$ or equivalently $A^d x=0$ for all $x\in\rd$ which contradicts the fact that $\rho(A)>0$. We conclude that $\lmax{(A^d)^\intercal R_I A^d}>0$ and so $\lmax{R_I}>1$. The same arguments combined with the inequality $\sum_{j=0}^{d-1} \left(A^\intercal\right)^j A^j\succeq \Idd+A^\intercal A$ can be used to prove $\lmax{\sum_{j=0}^{d-1} \left(A^\intercal\right)^j A^j}>1$. The inequalities $\alphaf(\beta)<\beta$ and $\betaf(\alpha)>\alpha$ follow readily from the first sentence. Let $\alpha,\beta>0$ and let $\gamma := \lmax{\sum_{j=0}^{d-1} \left(A^\intercal\right)^j A^j} \lmax{R_I}$. We have $\alpha\leq \alphaf(\beta)\iff \alpha\leq \beta/\gamma\iff \gamma\alpha\leq \beta\iff \betaf(\alpha)\leq \beta$.
\end{proof} 
The matrix norms are equivalent and to prove the compactness of $\sdp{\alpha,\beta}$ and the continuity of $G_k$, we can choose any matrix norm. We choose the Frobenius norm i.e. $B\mapsto \norm{B}_F:=\sqrt{\operatorname{tr}(B^* B)}$ which is submultiplicative, monotone ($B\succeq C$ implies that $\norm{B}_F \leq \norm{C}_F$) and unitarily invariant ($\norm{UBV}_F=\norm{B}_F$ for all unitary matrices $U,V$).
\begin{proposition}
\label{sdplemma}
The following assertions are valid:
\begin{enumerate}
\item Let $\alpha,\beta>0$. For all $\alpha'\leq \alpha$ and all $\beta\leq \beta'$: $\sdp{\alpha,\beta}\subseteq \sdp{\alpha',\beta'}$;
\item For all $\alpha>0$ and $\beta>0$, $\sdp{\alpha,\beta}$ is compact;
\item For all $\alpha>0$, $\sdp{\alpha,\betaf(\alpha)}$ is nonempty;
\item For all $\beta>0$, $\sdp{\alphaf(\beta),\beta}$ is nonempty;
\item Let $\alpha>0$ and $\alpha'\in(0,\alpha)$. Then :
\[
\{P\in\hp : \alpha' \Idd \preceq P-A^\intercal P A\preceq \alpha\Idd\} \subseteq \sdp{\alpha',\betaf(\alpha)}
\]
\item For all $\beta>0$:
\[
\lyap{A}\cap \lambda_{\rm max}^{-1}\left((0,\beta]\right)=\bigcup_{0<\alpha\leq \alphaf(\beta)} \sdp{\alpha,\beta}
\]
\end{enumerate}
\end{proposition}
\begin{proof}
\noindent ({\itshape 1}) Let $\alpha,\beta>0$, $\alpha'\leq \alpha$ and $\beta\leq \beta'$. The result is evident if $\sdp{\alpha,\beta}=\emptyset$. Suppose that $\sdp{\alpha,\beta}\neq\emptyset$ and let $P\in\sdp{\alpha,\beta}$. Then $P-A^\intercal P A\succeq \alpha \Idd\succeq \alpha'\Idd$ and $P\preceq \beta\Idd\preceq \beta' \Idd$. Finally, $P\in\sdp{\alpha',\beta'}$. 

\noindent ({\itshape 2}) Let $\alpha,\beta>0$. If $\setsdp=\emptyset$, the result obviously holds. Suppose that $\setsdp\neq \emptyset$ and let $P\in\setsdp$. We have $P\succeq 0$. This implies that $A^\intercal P A\succeq 0$ and thus $P\succeq P-A^\intercal P A$. This leads to $P\succeq \alpha \Idd$. Combining with $P\preceq \beta \Idd$, we conclude, since $\norm{\cdot}_F$ is monotone, that $\alpha \norm{\Idd}_F\leq \norm{P}_F \leq \beta \norm{\Idd}_F$ and $\setsdp$ is bounded. If a sequence $(P_n)_{n\in\nn}$ in $\setsdp$ converges to some $P$, then for all $x\neq 0$, $x^\intercal (P_n-A^\intercal P_n A)x$ converges to  $x^\intercal (P-A^\intercal P A)x$ and $x^\intercal P_n x$ converges to $x^\intercal Px$. We conclude that $P\in\setsdp$. Finally, $\setsdp$ is compact.

\noindent ({\itshape 3})-({\itshape 4}) By definition of $\alphaf$ and $\betaf$, and from Lemma~\ref{mori}, if we fix $\alpha>0$, then the unique $P\succ 0$ such that $P-A^\intercal P A=\alpha \Idd$ belongs to $\sdp{\alpha,\betaf(\alpha)}$. If we fix $\beta>0$ then the unique $P\succ 0$ such that $P-A^\intercal P A=\alphaf(\beta) \Idd$ belongs to $\sdp{\alphaf(\beta),\beta}$.

\noindent ({\itshape 5}) Let $\alpha>0$ and let $\alpha'<\alpha$. Let $P$ s.t $\alpha' \Idd \preceq P-A^\intercal P A\preceq \alpha \Idd$, then there exists $R$ s.t. $\alpha' \Idd \preceq R\preceq \alpha \Idd$ and $P-A^\intercal P A=R$. Hence $P=\sum_{k=0}^{+\infty} (A^\intercal)^k R A^k\preceq P_\alpha$ defined as the unique $P_\alpha\succ 0$ such that $P_\alpha-A^\intercal P_\alpha A=\alpha \Idd$. By definition of $\betaf(\alpha)$, $P\preceq  P_\alpha\preceq \betaf(\alpha)\Idd$. We conclude that $P\preceq \betaf(\alpha)\Idd$.

\noindent ({\itshape 6}) Let $\beta>0$ and $\alpha<\alphaf(\beta)$. Then $\sdp{\alpha,\beta}\subseteq  \lambda_{\rm max}^{-1}\left((0,\beta]\right)$. Moreover, as $P\in\sdp{\alpha,\beta}$ implies that $P-A^\intercal P A\succeq \alpha\Idd\succ 0$, we have $P\in\lyap{A}$. Finally, $\sdp{\alpha,\beta}\subseteq  \lambda_{\rm max}^{-1}\left((0,\beta]\right)\cap \lyap{A}$. Now, let $P\in\lambda_{\rm max}^{-1}\left((0,\beta]\right)\cap \lyap{A}$. There exists $R\succ 0$ s.t. 
 $P-A^\intercal P A=R\succeq \lmin{R}\Idd$. Then $P\in\sdp{\lmin{R},\beta}$. If $\lmin{R}\leq \alphaf(\beta)$, the result holds. Otherwise, the first statement of the proposition proves that for all $\alpha\leq \alphaf(\beta)<\lmin{R}$, $P\in\sdp{\alpha,\beta}$. This ends the proof. 
\end{proof}
\begin{lemma}
\label{lemcont}
The following statements hold:
\begin{enumerate}
\item  Let $B$ be in $\hddp$. For all $H\in\hdd$ such that $\norm{B^{-1}H}_F<1$ and $B+H\in\hddp$:
\[
\norm{(B+H)^{-1}-B^{-1}}_F\leq \sqrt{d}\dfrac{\norm{H}_F\norm{B^{-1}}_F}{1-\norm{B^{-1}H}_F} 
\]
\item The product mapping ${\hddp}^2\ni(B,C)\mapsto BC$ is continuous;
\item For all $B,C\in\hddp$, $|\lmax{B}-\lmax{C}|\leq \norm{B-C}_F$.
\end{enumerate}
\end{lemma}
\begin{proof}
The first statement comes from~\cite[Th 2.3.4]{golub2013matrix} which holds for $\norm{B}_2=\sqrt{\rho(B^* B)}$. Our extension follows from the fact that $\norm{B}_2\leq \norm{B}_F\leq \sqrt{d}\norm{B}_2$ for all $B\in\cdd$. The continuity of product mapping follows directly from the submultiplicative property of the Frobenius norm. For the last statement,~\cite[Corollary III 2.6]{bhatia1997matrix} asserts a stronger inequality with $\norm{\cdot}_2$ instead of $\norm{\cdot}_F$. Our conclusion follows from the fact that $\norm{\cdot}_F$ is greater than $\norm{\cdot}_2$.
\end{proof}
Lemma~\ref{lemcont} proves that the inverse, the matrix product and the maximal eigenvalue functions are continuous on $\hddp$. Thus their compositions are continuous.
\begin{proposition}
The function $G_k$ is continuous on $\lyap{A}$.
\end{proposition}
\begin{proof}
The function $G_k$ is the fraction of the logarithm of $P\mapsto H_{P,\topt{P}}^{-1}(\nu_k)$ and the one of $P\mapsto \norm{A}_P$. Then it is continuous if $P\mapsto H_{P,\topt{P}}^{-1}(\nu_k)$ and $P\mapsto \norm{A}_P$ are. So, to prove that $G_k$ is continuous, we have to show that $P\mapsto \mu(P)$, $P\mapsto \norm{q}_P$, $P\mapsto \topt{P}$ and $P\mapsto \norm{A}_P$ are continuous. 

First, $\mu$ is continuous since the mappings $P\mapsto x^\intercal P x$ are linear for all $x\in\xin$. For all $P\succ 0$, $x\mapsto x^\intercal P x$ are convex and $\mu(P)$ can be expressed as $\max_{x\in \mathcal{E}(\xin)} x^\intercal P x$ where $\mathcal{E}(\xin)$ denotes the finite set of vertices of $\xin$. Finally, $\mu$ is continuous as the maximum of a finite family of continuous functions. When $q=0$, $\norm{q}_P^*$ is obviously continuous. When $q\neq 0$, Lemma~\ref{lemcont} proves that inverse is continuous and thus $P \mapsto \norm{q}_P^*=\sqrt{q^\intercal P^{-1} q}$ is continuous on $\hddp$ as the composition of continuous functions. Using again Lemma~\ref{lemcont}, the functions $P\mapsto \topt{P}=\lmax{Q P^{-1}}$ and $P\mapsto \norm{A}_P=\lmax{A^\intercal PA P^{-1}}$ are continuous by composition. 
\end{proof}
\begin{theorem}
\label{thsdpopt}
Let $k\in\pos_\nu$. We define, for all strictly positive scalars $\alpha$ and $\beta$ the nonlinear semidefinite program:
\begin{equation}
\label{SDPQgen}\tag{${\rm SDP}_{\alpha,\beta}$}
\Min\  G_k(P)\ \st \ P\in\setsdp
\end{equation}
Let $\alpha,\beta>0$ be such that Problem~\eqref{SDPQgen} is feasible. Then :
\begin{enumerate}
\item Problem~\eqref{SDPQgen} admits an optimal solution $P^{\alpha,\beta}$;
\item For all $0<\varepsilon<\varepsilon_k(P^{\alpha,\beta})$, all $\varepsilon$-optimal solutions $P_\varepsilon^{\alpha,\beta}$ satisfy:
\[
\klyap_\nu\left(k,P_\varepsilon^{\alpha,\beta},\topt{P_\varepsilon^{\alpha,\beta}}\right)=\klyap_\nu\left(k,P^{\alpha,\beta},\topt{P^{\alpha,\beta}}\right)
\] 
\end{enumerate}
Moreover, for all $\beta>0$, there exists $\alpha\in (0,\alphaf(\beta)]$ such that:
\[ 
\klyap_\nu\left(k,P^{\alpha,\beta},\topt{P^{\alpha,\beta}}\right)=\min_{P\in\lyap{A}} \klyap_\nu\left(k,P,\topt{P}\right)
 \]
\end{theorem}
\begin{proof}
\noindent ({\itshape 1}) The function $G_k$ is continuous on the compact set $\setsdp$ and thus attains its minimum on it.

\noindent ({\itshape 2}) This is a direct application of Prop.~\ref{optimerror} with $\overline{P}=P^{\alpha,\beta}$.
Let $\beta$ be a positive scalar. The function $\lambda_{\rm max}$ is an homogeneous function s.t. $\lmax{\hddp}\subseteq \rr_+^*$ and then we have $\lyopt_k=V_k^{\lambda_{\rm max},\beta}$ which is finite. From Prop.~\ref{propopt}, we can find $\overline{P}\in\lyap{A}$ such that $\lfloor G_k(\overline{P})\rfloor=\lfloor\lyopt_k\rfloor=\lfloor V_k^{\lambda_{\rm max},\beta}\rfloor$. Let $\varepsilon<\varepsilon_k(\overline{P})$. Hence there exists $P_{\varepsilon}\in \lyap{A}\cap \lambda_{\rm max}^{-1}\left((0,\beta]\right)$ s.t.
$V_k^{\lambda_{\rm max},\beta}\leq G_k(P_{\varepsilon})\leq V_k^{\lambda_{\rm max},\beta}+\varepsilon$ and $\lfloor G_k(P_{\varepsilon})\rfloor)=\lfloor V_k^{\lambda_{\rm max},\beta}\rfloor$. From the last statement of Lemma~\ref{sdplemma}, there exists $0<\alpha\leq \alphaf(\beta)$ such that $P_{\varepsilon}\in \sdp{\alpha,\beta}$. This implies that
$V_k^{\lambda_{\rm max},\beta}\leq G_k(P^{\alpha,\beta})\leq G_k(P_{\varepsilon})\leq V_k^{\lambda_{\rm max},\beta}+\varepsilon$
and $P_{\varepsilon}$ is an $\varepsilon$-solution with $\varepsilon<\varepsilon_k(\overline{P})$ and we conclude from the second statement that $\lfloor G_k(P^{\alpha,\beta})\rfloor=\lfloor G_k(P_{\varepsilon})\rfloor=\lfloor V_k^{\lambda_{\rm max},\beta}\rfloor$.
\end{proof}
\begin{corollary}
\label{finalcoro}
For all $\beta>0$, $\displaystyle{\lyopt_k=\lim_{\alpha\downarrow 0} G_k(P^{\alpha,\beta})}$
\end{corollary}
Theorem~\ref{thsdpopt} and Corollary~\ref{finalcoro} prove that we can choose any $\beta$ we want large enough to avoid numerical issues in an implementation perspective. The problem is to guess a value $\alpha>0$ small enough to obtain the optimal value $\lfloor \lyopt_k\rfloor$. We could design a kind of descent algorithm using the fifth statement of Prop.~\ref{sdplemma} to produce a strict decreasing sequence $(G_k(P^{\alpha_n,\beta}))_{n\in\nn}$. The main issue is to construct the stopping condition. This relies on a first order condition. 

This theoretical approach proves that we can obtain a matrix $P\in\lyap{A}$ that achieves $\lyopt_k$ and is the optimal solution of a nonlinear semidefinite program for which the constraints set is a compact subset of $\lyap{A}$. The problem is to solve the associated nonlinear semidefinite program. The main difficulty comes from the objective function. The objective function can be nonconvex (see Figure~\ref{nonconvex}) and nondifferentiable ($\lmax{\cdot}$ and $\mu$ are not Fréchet differentiable). The first difficulty can be overcame using recent methods to solve nonconvex and general nonlinear semidefinite programs such as primal-dual type methods~\cite{andreani2020optimality} or sequencial quadratic programming approaches~\cite{yamakawa2022stabilized}. The interesting reader could also consult the survey paper~\cite{Yamashita2015}. However, those methods require sufficiently smooth objective function and need its gradient. To develop an extension of those methods specific to our problem is beyond the scope of this paper and left to future works. 
\subsection{In Practice}
In practice, we adopt a more naive approach based on the classical linear semidefinite programming. We compute a matrix $P$ in $\lyap{A}$ as an optimal solution of a linear semidefinite program. Recall that $\mathcal{E}(\xin)$ denotes the set of vertices of $\xin$. We tried five objective functions :
\begin{equation}
\label{objsdp}
\begin{array}{c}
F_0(P):=\displaystyle{\max_{x\in\mathcal{E}(\xin)} x^\intercal P x};\  F_1(P):=F_0(P-Q);\
F_2(P):=\displaystyle{\sum_{x\in\mathcal{E}(\xin)} x^\intercal P x};\\ F_3(P):=\lmax{P} \text{ and } F_4(P):=0
\end{array}
\end{equation} 
Obviously, the fifth objective function leads to a feasibility semidefinite problem. For each objective functions, we add the Lyapunov condition up to perturbation $\varepsilon=0.01$ to ensure the strict positiveness of discrete Lyapunov equation matrix constraint, we get the following family of problems for $i=0,1,2,3,4$:
\begin{equation}
\label{eq:klmaxp}
\Min  F_i(P)\ \st\ 
         P-A^\intercal P A-\varepsilon \Idd\succeq 0\text{ and } P\succeq 0
\end{equation}
Since the constraints set of Pb~\eqref{eq:klmaxp}, they have an optimal solutions when objective functions are coercive. This is the case for $i=3$ and the case for $i=0,1,2$ when $\xin$ is full-dimensional. Problems~\eqref{eq:klmaxp} with $i\in\{0,1\}$ can be rewritten as linear semi-definite programs. Indeed, to minimize $\max_{x\in\mathcal{E}(\xin)} x^\intercal B x$ is equivalent to minimize $t$ such that $t\geq x^\intercal B x$ for all $x\in\mathcal{E}( \xin)$. For each $i=0,1,2,3,4$, we will write $\pmaxi$ for an optimal solution of Pb.~\eqref{eq:klmaxp}. Once $\pmaxi$ computed, we have to evaluate $\topt{\pmaxi}$, $\mu(\pmaxi)$, $\norm{q}_{\pmaxi}^*$ and $\norm{A}_{\pmaxi}$. To evaluate $\topt{\pmaxi}$ boils down to compute the maximal eigenvalue of $Q P_{\rm max}^{{i}^{-1}}$ which can be done by any solver of linear algebra. It can be also computed from semidefinite programming solvers. To evaluate $\mu(\pmaxi)$ can be done by the comparison of values of $x\mapsto x^\intercal P x$ on each vertex $x$ of $\xin$ again, this operation can be easily implemented. In a future work, we consider scalable techniques (e.g.~\cite{konno1976maximization}) to compute the exact value of $\mu(P)$. The value $\norm{q}_{\pmaxi}^*$ is the square root of $q^\intercal P_{\rm max}^{{i}^{-1}} q$ and thus can be done using any linear algebra solver. Finally, the value $\norm{A}_{\pmaxi}$ can be computed similarly to $\topt{\pmaxi}$.
\section{Comparison with the Diagonalizable Case}
\label{comparison}
In this part, we compare the Lyapunov approach with the case where the matrix $A$ is diagonalizable. We will prove that the formula proposed in this paper generalizes the one constructed in the diagonalizable case and recalled at Eq.~\eqref{oldformula}. In~\cite{DBLP:journals/jota/Adje21}, the dependency with respect to the eigenvectors matrix $U$ was not explicit, here, we make it explicit. Let us consider :
\[
\eigen(A):=\left\{U\in\cdd \text{ non-singular } : \exists\, D\in\cdd \text{ diagonal }\text{ s.t } A=UDU^{-1}\right\} 
\]
We suppose in this part that $\eigen(A)\neq \emptyset$. 
\begin{lemma}
\label{lemmaequal}
Let $U\in\eigen(A)$. Then:
\begin{enumerate}
\item $\lmax{\left((UU^*)^{-1}\right)^{-1/2}Q\left((UU^*)^{-1}\right)^{-1/2}}=\lmax{U^*QU}$;
\item $(UU^*)^{-1}\in\lyap{A}$;
\item $\norm{A}_{(UU^*)^{-1}}=\rho(A)$;
\item $\norm{q}_{(UU^*)^{-1}}^*=\norm{U^*q}_2$.
\end{enumerate}
\end{lemma}

\begin{proof}
\label{uulyap}
To simplify the notations, we will write $U^{-*}$ for ${U^{*}}^{-1}$.

\noindent {\itshape 1.} We have $\lmax{\left((UU^*)^{-1}\right)^{-1/2}Q\left((UU^*)^{-1}\right)^{-1/2}}=\lmax{Q\left((UU^*)^{-1}\right)^{-1/2}\left((UU^*)^{-1}\right)^{-1/2}}=\lmax{QUU^*}=\lmax{U^*QU}$.

\noindent {\itshape 2.}  Let $W:=(UU^*)^{-1}$. Clearly $W\in\hddp$ and for some diagonal matrix $D$, we have $W-A^\intercal W A=W-A^* WA=W-U^{-*}D^*U^*U^{-*}U^{-1}UDU^{-1}=W-U^{-*}D^*DU^{-1}=U^{-*}(\Idd-D^*D)U^{-1}$. As $\rho(A)<1$ then $(\Idd-D^*D)\in\hddp$ and $W\in\lyap{A}$.

\noindent {\itshape 3.} To prove that $\norm{A}_{(UU^*)^{-1}}=\rho(A)$, we apply a  change of variable by setting $y=U^{-1}x$ and, for some diagonal matrix $D$, we obtain $x^*A^* (UU^*)^{-1} Ax=y^*U^*U^{-*}D^*U^*U^{-*}U^{-1}UDU^{-1}Uy=D^*D$ and :
\[
\displaystyle{\norm{A}_{(UU^*)^{-1}}^2}=\displaystyle{\max_{\substack{x\in\cd\\ x\neq 0}}  \dfrac{x^* A^* (UU^*)^{-1} A x}{x^*(UU^*)^{-1}x}}=\displaystyle{\max_{\substack{y\in\cd\\ y\neq 0}}  \dfrac{y^*D^*Dy}{y^* y}}=\lmax{D^* D}=\rho(A)^2
\]

\noindent {\itshape 4.} We have $\displaystyle{\norm{q}_{(UU^*)^{-1}}^*=\max_{\norm{y}_{(UU^*)^{-1}=1}} |q^* y|=\max_{\norm{z}_2=1} |q^* U z|=\norm{U^* q}_2}$.
\end{proof}
Lemma~\ref{lemmaequal} shows that, for all $U\in\eigen(A)$, $(UU^*)^{-1}\in\lyap{A}$ and $\topt{(UU^*)^{-1}}=\max\{\lmax{U^*QU},0\}$. We can define the homeomorphism $H_{(UU^*)^{-1},\topt{(UU^*)^{-1}}}$ following Eq~\eqref{homeoseq}. We suppose that Assumption~\eqref{ineqforall} still holds. This implies that we must have $\nu_0<  H_{(UU^*)^{-1},\topt{(UU^*)^{-1}}}(1)$. 

We can improve $\kdiag_\nu$ formula to refine the case where $Q\preceq 0$ as it is done in this paper (the case $t=0$). Indeed, at Prop.~\ref{nukineg}, the quadratic part $x^\intercal {A^\intercal }^k Q A^k x$ can be upper-bounded by 0 when $Q\preceq 0$ by setting $t=0$ whereas in~\cite{DBLP:journals/jota/Adje21} it was upper-bounded by $|\lmax{U^* Q U}| \norm{A^k x}_2^2$. We define, for $U\in\eigen(A)$:
\begin{equation}
\label{knew}
\kdiag_{{\rm new},\nu}(k,U)=
 \left\lfloor\dfrac{1}{\ln(\rho(A))}\ln\left(H_{(UU^*)^{-1},\topt{(UU^*)^{-1}}}^{-1}(\nu_k)\right)\right\rfloor+1
\end{equation}
Let $U$ be in $\eigen(A)$ and $k$ be in $\pos_\nu$. From Lemma~\ref{lemmaprmin}, we also have :
\[
\kdiag_{{\rm new},\nu}(k,U)=\left\lfloor \dfrac{1}{\ln(\rho(A))}\ln\left( h_{(UU^*)^{-1},\topt{(UU^*)^{-1}}}^-(\nu_k)\right)\right\rfloor+1;
\]

\begin{theorem}
Let $U\in\eigen(A)$, we have, for all $k\in\pos_\nu$:
\[\klyap_\nu(k,(UU^*)^{-1},\topt{(UU^{*})^{-1}})=\kdiag_{{\rm new},\nu}(k,U)\leq \kdiag_\nu(k,U)\]
\end{theorem}
\begin{proof}
The equality is a consequence of Lemma~\ref{lemmaequal}. The inequality follows directly from Lemmas~\ref{propcont} and~\ref{lemmaprmin}. Indeed, $|\lmax{U^*QU}|$ belongs to $\cont((UU^*)^{-1})$ and $\topt{(UU^*)^{-1}}=\max\{\lmax{U^*QU},0\}\leq |\lmax{U^*QU}|$.
Moreover, simple reformulations of  Eq~\eqref{oldformula}  permit to prove that 
$h_{(UU^*)^{-},|\lmax{U^*QU}|}^-(\nu_k)=\left(\sqrt{\nu_k+\vdiag^2}-\vdiag\right)\left(\sqrt{|\lmax{U^*QU}|\mus{(UU^{*})^{-1}}}\right)^{-1}$.
\end{proof}

\begin{corollary}[Formula improvment]
Let $U\in\eigen(A)$ and $k\in\pos_\nu$ . The following statement holds:
\begin{enumerate}
\item $\kdiag_{{\rm new},\nu}(k,U)$ is well-defined i.e a natural integer;
\item $\begin{array}{ll}\kdiag_{{\rm new},\nu}(k,U)&=\min\left\{j\in\nn : H_{(UU^*)^{-1},\max\{\lmax{U^* QU},0\}}(\rho(A)^j)< \nu_k\right\}\\
&=\min\left\{j\in\nn : \rho(A)^j< H_{(UU^*)^{-1},\max\{\lmax{U^* QU},0\}}^{-1}(\nu_k)\right\}\end{array}$;
\item $\kdiag_{{\rm new},\nu}(k,U)\geq k+1$;
\item For all $j\geq \kdiag_{{\rm new},\nu}(k,U)$, $\nu_j< \nu_k$;
\item Let $j\in\pos_\nu$. If $\nu_j\leq \nu_k$ then $\kdiag_{{\rm new},\nu}(k,U)\leq \kdiag_{{\rm new},\nu}(j,U)$.
\end{enumerate}
\end{corollary}
The results of Prop.~\ref{inter} can be also extended to $\kdiag_{{\rm new},\nu}(k,U)$. 
\section{Examples}
\label{benchexamples}
\subsection{Running Example from~\cite{DBLP:journals/jota/Adje21}}
We consider the example presented in~\cite{DBLP:journals/jota/Adje21}. It concerns the discretisation of an harmonic oscillator $\ddot{x}+\dot{x}+x=0$ by an explicit Euler scheme. The discretization step is set to 0.01. We introduce the position variable, $x$ and the speed variable $v$. We assume that the initial conditions can be taken into the set $[-1,1]^2$. The Euler scheme becomes a linear discrete-time system in dimension two defined as follows: 
\begin{equation}
\label{harmonic}
\begin{pmatrix}
x_{k+1}\\
v_{k+1}
\end{pmatrix}
=\begin{pmatrix}
1 & 0.01\\
-0.01 & 0.99
\end{pmatrix}
\begin{pmatrix}
x_{k}\\
v_{k}
\end{pmatrix},\ (x_0,v_0)\in [-1,1]^2
\end{equation}
For this linear system, we are interested in computing:
\begin{itemize}
\item the maximal value of the Euclidean norm of the state-variable $\norm{(x_k,v_k)}_2^2$ and thus $Q=\begin{pmatrix}
1 & 0\\ 
0 & 1
\end{pmatrix}$ ;
\item the square of the position variable $x_{k}^2$ and thus $Q=\begin{pmatrix}
1 & 0\\ 
0 & 0
\end{pmatrix}$;
\item the speed variable $v_{k}^2$ and thus $Q=\begin{pmatrix}
0 & 0\\ 
0 & 1
\end{pmatrix}$.
\end{itemize}
In~\cite{DBLP:journals/jota/Adje21} where the method was developed for the diagonalizable case, we used a matrix $U\in\eigen(A)$ without any optimization procedure. The matrix $U$ is composed of the eigenvectors with one on the first coordinate. From this matrix $U$ we were able to compute $\kdiag_\nu(k,U)$. Adapting the second statement of Prop.~\ref{inter} to $\kdiag_\nu(k,U)$, its smallest value is achieved at $k=\bigkse_\nu$. We recall the results for the three later matrices :
\begin{itemize}
\item For $Q=\begin{psmallmatrix}
1 & 0\\ 
0 & 1
\end{psmallmatrix}$, we had $\kopt=0$, $\nuopt=2$ and $\kdiag_\nu(0,U)=111$.
\item For $Q=\begin{psmallmatrix}
1 & 0\\ 
0 & 0
\end{psmallmatrix}$, we had $\kopt=61$, $\nuopt\simeq 1.64886$ and $\kdiag_\nu(61,U)=90$.
\item For $Q=\begin{psmallmatrix}
0 & 0\\ 
0 & 1
\end{psmallmatrix}$, we had $\kopt=0$, $\nuopt=1$ and $\kdiag_\nu(0,U)=140$.
\end{itemize}
For each matrix $Q$, we compute solutions $\pmaxi$ of semidefinite Problems~\eqref{eq:klmaxp}. Those matrices are computed once since their computation does not depend on $k$. Then at each $k$, we compare their associated value $\klyap_\nu(k,\pmaxi,\topt{\pmaxi})$ and pick the smallest.

Now we present the computations of $\klyap_\nu(k,P,t)$ for the three latter matrices using the proposed method of the paper. For $Q=\begin{psmallmatrix}
1 & 0 \\ 0  &1
\end{psmallmatrix}$, we compute a solution of the discrete Lyapunov equation only for $k=0$ since for all $k\in\nn$, $\nuopt=\nu_0\geq \nu_k$. The comparison between the five possibilities proposed by the optimal solutions of the semidefinite programs~\eqref{eq:klmaxp} leads to pick:
\[
\bP=
\begin{pmatrix}
   4.9501 +      0i &  1.5475 + 0.0000i\\
   1.5475 - 0.0000i  & 4.5313 +      0i
\end{pmatrix}
\]
which corresponds to the matrix obtained by solving the feasibility problem i.e. $i=4$. From this matrix $\bP$, we then compute $\topt{\bP}$:
\[
\topt{\bP}=\lmax{Q\bP^{-1}}\simeq 0.31456
\]
and the maximal value $\mu(\bP)$ of $x\mapsto x^\intercal \bP x$ on $\xin$:
\[
\mu(\bP)
=\displaystyle{\max\left\{\begin{pmatrix}1\\1\end{pmatrix}^\intercal \bP\begin{pmatrix}1\\1\end{pmatrix},\begin{pmatrix}-1\\-1\end{pmatrix}^\intercal \bP\begin{pmatrix}-1\\-1\end{pmatrix},\begin{pmatrix}-1\\1\end{pmatrix}^\intercal \bP\begin{pmatrix}-1\\1\end{pmatrix},\begin{pmatrix}1\\-1\end{pmatrix}^\intercal \bP\begin{pmatrix}1\\-1\end{pmatrix}\right\}}
\simeq 12.57646
\]
and finally, the norm of $A$ with respect to $\bP$:
\[
\norm{A}_{\bP}=\sqrt{\lmax{A^\intercal \bP A  \bP^{-1}}}\simeq 0.996922
\]
Note that the spectral radius of $A$ is about $0.99504$. Finally, this leads to :
\[
\klyap_\nu(0,\bP,\topt{\bP})
=\left\lfloor\dfrac{1}{\ln(\norm{A}_{\bP})}\dfrac{\sqrt{2}}{\sqrt{\topt{\bP}}\mu(\bP)}\right\rfloor+1=\lfloor 110.63173\rfloor+1=111
\]
The value is the one we found with the previous method proposed in~\cite{DBLP:journals/jota/Adje21}.

Now we consider $Q=\begin{psmallmatrix}
1 & 0 \\ 0  &0
\end{psmallmatrix}$. We pick the same matrix for all $k$ since it still provides the minimal value for $P\mapsto \klyap(k,P,\topt{P})$ with respect to five possible matrices. Again this corresponds to the matrix obtained by solving the feasible problem i.e. $i=4$. We recover the matrix $\bP$ given above used for $Q=\begin{psmallmatrix}
1 & 0 \\ 0  &1
\end{psmallmatrix}$.
Compared to the case where $Q=\begin{psmallmatrix}
1 & 0 \\ 0  &1
\end{psmallmatrix}$, the only change is in the computation of $\topt{\bP}=\lmax{Q\bP^{-1}}$. For $Q=\begin{psmallmatrix}
1 & 0 \\ 0  &0
\end{psmallmatrix}$, we have $\lmax{Q\bP^{-1}}\simeq 0.226165$.
The values $\mu(\bP)$ and $\norm{A}_P$ does not depend on $Q$ and are the one computed previously. Some values of 
$\klyap_\nu(k,\bP,\topt{\bP})$ are given in Table~\ref{oldrun}.
\begin{center}
\begin{table}[h!]
{\footnotesize
\begin{center}
\begin{tabular}{|*{8}{c|}}
\hline
$k$&  0 & 1& 10 & 20 & 30 & 50 & {\bf 61} \\
\hline
$\nu_k$ & 1 & 1.0201 & 1.19055& 1.353 & 1.4818  &1.62944 & {\bf 1.64886}\\
\hline
$\klyap_\nu(k,\pmax,t)$&  170 & 167  & 142 & 121 & 106  &91 & {\bf 89}\\
\hline
\end{tabular}
\end{center}
}
\caption{Some values of $\klyap_\nu(k,\bP,\topt{\bP})$ for the running example of~\cite{DBLP:journals/jota/Adje21}}
\label{oldrun}
\end{table}
\end{center}
We have gained one iteration with respect to~\cite{DBLP:journals/jota/Adje21}.

Finally, we consider $Q=\begin{psmallmatrix}
0 & 0 \\ 0  &1
\end{psmallmatrix}$. We compute a solution of the discrete Lyapunov equation only for $k=0$ since for all $k\in\nn$, $\nuopt=\nu_0\geq \nu_k$. We use the same matrix as before; the feasibility problem provides again the smallest integer at $k=0$. Again, the only change is in the computation of $\topt{\bP}=\lmax{Q\bP^{-1}}\simeq 0.247068$ which gives $\klyap_\nu(0,\bP,\topt{\bP})=184$ which is worse than the one computed from the diagonalizable case.
\subsection{Higher Random Examples}
In this subsection, we propose to test our technique on randomly generated systems. Those benchmark tests are done using Julia 1.7.3~\cite{Julia-2017} for which our problems were more quickly solved than Octave 6.3.0. For those tests, we construct a specific family of initial polytope sets parameterized by the dimension of the system. The other data such as the system matrix and the objective function are randomly generated.

\subsubsection{Initial Polytope $\xin$ Generation}
The proposed method needs the existence of a strictly positive term $\nu_k$. To ensure its existence, we extend Prop. 3.5 of~\cite{DBLP:journals/jota/Adje21}.

\begin{proposition}[Extension of Prop. 3.5 of~\cite{DBLP:journals/jota/Adje21}]
Suppose that $0\in\inte\xin$. If $Q\succeq 0$;
or if $q\neq 0$ and $Q\succeq -\norm{q}_2^{-1} qq^\intercal$, then $\displaystyle{\nu_0=\max_{x\in\xin} x^\intercal Q x+q^\intercal x>0}$.
\end{proposition}
\begin{proof}
The case $Q\succ 0$ is treated in~\cite[Prop. 3.5]{DBLP:journals/jota/Adje21}. Assume that $q\neq 0$ and $Q\succeq -\norm{q}_2^{-1} qq^\intercal$. As $0\in\inte\xin$, there exists $\varepsilon\in (0,1)$ such that $z_\varepsilon:=\varepsilon q\norm{q}_2^{-1}\in\xin$. Finally, $\nu_0 \geq z_\varepsilon^\intercal Q z_\varepsilon +q^\intercal z_\varepsilon \geq -\varepsilon^2\norm{q}_2^{-3} (q^\intercal q q^\intercal q) +\varepsilon \norm{q}_2^{-1}q^\intercal q= (1- \varepsilon)\norm{q}_2>0$.
\end{proof}

Now, it suffices to construct an initial polytope with 0 in its interior and with a small number of vertices. Any box of the form $\Pi_{i=1}^d [-a_i,a_i]$ with $a_i>0$ contains 0 in its interior. However, the number of vertices grows exponentially with respect to the dimension of the space. The simplest situation it to use, as an initial polytope, a $d$-simplex i.e. the convex hull $d+1$ affinely independant vectors. Indeed, as in $\rd$, a convex set has nonempty interior if and only if its dimension is equal to $d$ (e.g see Prop.11.2.27~\cite{berger}), we need at least $d+1$ affinely independant vectors for our use. Then we define the family $\{x_k\}_{k=1,\ldots,d+1}$ of vectors of $\rd$ where :
\[
\begin{array}{lr}
x_{1,i}= -1, \forall\, i=1,\ldots d,& x_{2,i}=\left\{\begin{array}{cr} 1 & \text{ if } i=1\\ 0 & \text{ if } 2\leq i\leq d\end{array}\right.,\\
\\
 x_{k,i}=\left\{\begin{array}{cr} -1 & \text{ if } 1\leq i\leq k-2\\ 1 & \text{if } i=k-1 \\ 0 & \text{ if } k\leq i\leq d\end{array}\right. \text{ and } &  x_{d+1,i}=\left\{\begin{array}{cr} -1 & \text{ if } 1\leq i\leq d-1\\ 1 & \text{if } i=d \end{array}\right.
\end{array}
\]
and we will use :
\begin{equation}
\label{initeq}
\xin=\conv(\{x_k,k=1,\ldots,d+1\})
\end{equation}
\begin{proposition}
\label{initbench}
We have:
\begin{enumerate}
\item The family $\{x_k\}_{k=1,\ldots,d+1}$ is affinely independant;
\item Let us define the matrix $F$ such that for all $i,j\in \{1,\ldots d\}$:
\[
F_{i,j}=\left\{
\begin{array}{lr} 
2^{j-i-1} & \text{ if } i>j\\
-\frac{1}{2} & \text{ if } i=j\\
0 & \text{ if }j > i
\end{array}
\right.
\]
Then :
\[
y\in\conv\left(\left\{x_k,k=1,\ldots,d+1\right\}\right)\iff \left\{\begin{array}{l}\displaystyle{F_i y\leq 2^{-i}},\ \forall\, i=1,\ldots,d\\ 
\\
\displaystyle{\sum_{j=1}^d} 2^{j-1} y_j\leq 1\end{array}\right.
\]
\item $0\in \inte\conv\left(\left\{x_k,k=1,\ldots,d+1\right\}\right)$.
\end{enumerate}
\end{proposition}
\subsubsection{Random Generation of Other Data}
Now to construct our test problems, we need :
\begin{enumerate}
\item a matrix $A$ with a spectral radius strictly smaller than 1;
\item a matrix $Q$ and a vector $q$ to define our quadratic objective function.
\end{enumerate}
To generate the matrix $A$, we adopt the following methodology:
we generate a random matrix $A$ and if its spectral radius $r_A$ is greater than one, then we replace $A$ by $A/(r_A+\varepsilon)$ where $\varepsilon=0.05, 0.1, 0.5, 1$ or $2$. The final produced matrix has a spectral radius strictly smaller than one. The different lines of Tables~\ref{convextable},~\ref{concavetable} and~\ref{lineartable} for each system dimension represent the different values for $\varepsilon$. The random generation for the vector $q$ is standard. For the matrix $Q$, it depends on the convexity that we want for the problem. When the objective function has to be convex, we generate $Q$ randomly and it the result is not definite positive, we replace $Q$ by $Q-\lambda_{\rm min}(Q)\Idd$. When the objective function has to be linear then $Q=0$ and when the objective function is concave $Q=-\norm{q}_2^{-1} qq^\intercal$ to have at least $\nu_0>0$. 
\subsubsection{Benchmarks Results}
We divide the benchmark tests in three main subclasses :
\begin{itemize}
\item the problems with a nonlinear convex quadratic objective function;
\item the problems with a nonlinear concave quadratic objective function;
\item the problems with a linear objective function.
\end{itemize}
For each subclass, we are interested in computations of : 
\begin{itemize}
\item $\kopt$, the smallest integer $k$ such that  $\nu_k=\max_{x\in\xin} x^\intercal {A^k}^\intercal  Q A^k x+q^\intercal A^k x$;
\item $\klyap_\nu$, an integer after which the search of the maximizer is useless and such that $\klyap_\nu-1$ is a safe number of quadratic problems to solve to compute $\nuopt$.
\end{itemize}
Those tests aim to analyze the behaviours of $\kopt$ and $\klyap_\nu$ with respect to the dimension of the system and the spectral radius of the matrix system. Therefore, for all dimensions $d\in\{3,5,10,20,30\}$ of the system and all possible pertubations of the spectral radius $\varepsilon\in\{0.05,0.1,0.5,1,2\}$, we generate 100 instances of the problem. The 100 instances differ from the system matrix and the objective function. The initial polytope is deterministic and depends on the dimension as it is shown at Eq.~\eqref{initeq}. Tables~\ref{convextable},~\ref{concavetable} and~\ref{lineartable} are divided in the following columns:
\begin{itemize}
\item $d$ : the dimension of the system;
\item $\varepsilon$ : the fixed deviation over the spectral radius of the randomly generated matrix system;
\item Spec A : Mi/Av/Mx = the minimum/average/maximum spectral radius of the 100 final generated system matrices;  ;
\item $\kopt$ : Mi/Av/Mx = the minimum/average/maximum integer $\kopt$ computed over the 100 final generated instances;
\item $\klyap_\nu$ : Mi/Av/Mx  = the minimum/average/maximal integer $\klyap_\nu$ computed over the 100 generated instances;
\item $\klyap_\nu-\kopt$ : Mi/Av/Mx = the minimum/average/maximum difference between $\klyap_\nu$ and $\kopt$ computed over the 100 final generated instances;
\end{itemize} 
The values in the column $\klyap_\nu-\kopt$ indicates how many supplementary useless quadratic maximizations we solve.

\renewcommand{\arraystretch}{1.5}

\begin{table}[h!]
\centering
{\footnotesize
\begin{tabular}{|>{\centering }p{0.25cm}|>{\centering }p{0.5cm}|>{\centering }p{0.7cm}>{\centering }p{0.7cm}>{\centering }p{0.7cm}|>{\centering }p{0.2cm}>{\centering }p{0.3cm}>{\centering }p{0.3cm}|>{\centering }p{0.4cm}>{\centering }p{0.6cm}>{\centering }p{0.7cm}|>{\centering }p{0.5cm}>{\centering }p{0.65cm}>{\centering\arraybackslash }p{0.7cm}|}
\hline
\hline
\hline
\multicolumn{14}{|c|}{Convex}\\
\hline
d & $\varepsilon$ & \multicolumn{3}{c|}{Spec A}& \multicolumn{3}{c|}{$\kopt$} & \multicolumn{3}{c|}{$\klyap_\nu$}  &\multicolumn{3}{c|}{$\klyap_\nu-\kopt$}\\
\hline
& &Mi & Av & Mx & Mi & Av & Mx & Mi & Av & Mx &  Mi & Av & Mx \\
\hline
\hline
\multirow{5}{*}{3} & 0.05 &0.9762 & 0.9892& 0.9943
&0 &2.2 &17 & 13& 183.9&1534 & 13&181.7 &1530 \\
\cline{2-14}
&0.1 &0.9483& 0.9787& 0.9883& 0& 1.6& 7& 7 & 132.1 &2715 &7 & 130.5&2709  \\
\cline{2-14}
&0.5 &0.785 & 0.9075&0.9512 & 0&0.7 &3
&1 &11.9 & 80& 1 &11.2 & 80\\
\cline{2-14}
&1 &0.6713& 0.8186& 0.8915& 0& 0.6& 4& 1 & 6.6&27 &1 & 6&26  \\
\cline{2-14}
&2 & 0.5329 & 0.7139 &0.8076 & 0&0.2 &2
&1 &2.7 & 8 & 1 &2.5 & 7\\
\hline
\hline
\multirow{5}{*}{5} &0.05 &0.9884 &0.9923 & 0.9947& 0&2.2 &18 &35 &415.5 & 3577& 34&413.2 &3559\\
\cline{2-14}
&0.1 &0.9727 &0.9855 &0.9921 &0 &1.7 &10 &22 & 241.1&8790 &22 &242.3 &8781 \\
\cline{2-14}
&0.5 & 0.8853& 0.9272&0.9536 & 0&0.9 &5 & 5& 26.4& 162& 5&25.4 &159\\
\cline{2-14}
&1 &0.7751&0.865&0.9105 &0 &0.7 &4 &2 & 15.6&181 &2&14.9 &181 \\
\cline{2-14}
&2 & 0.6582& 0.7729 &0.8598 & 0&0.2&4 & 1& 5.4& 22& 1&5.1 &20\\
\hline
\hline
\multirow{5}{*}{10}&0.05  & 0.9932 &0.9948 &0.9963 & 0&2.6 &39 &81 &1025.9 & 5959& 81& 1023.2&5948 \\
\cline{2-14}
&0.1 & 0.9866& 0.9895 & 0.9933&0 &3 &40 &16 &599.5 & 7149&16 &596.4 &7109 \\
\cline{2-14}
&0.5 & 0.9287& 0.9494& 0.9677&0 &1.4 &18 &5 &96.1 &2633 &5 &94.6 &2615 \\
\cline{2-14}
&1 & 0.8756&0.9027 &0.9271 &0 &1 &6 &2 &31.7 &174 &2 &30.6 &168 \\
\cline{2-14}
&2 & 0.7934& 0.827& 0.8619&0 &0.3 &6 &3 &11.2 &48 &3 &10.8 &42 \\
\hline
\hline
\multirow{5}{*}{20}&0.05 &0.9955 &0.9963 &0.9973 & 0&3.6 &36 &224 & 2240.9& 25832& 224&2237.2 &25805 \\
\cline{2-14}
&0.1 &0.9911 &0.9926 &0.9942 &0 &3.2 &27 &82 &1208.2 & 18159& 82& 1205& 18132\\
\cline{2-14}
&0.5 & 0.9565& 0.9647& 0.9709& 0&1.4 &7 &15 &111.8 &614 &15 &110.4 &607 \\
\cline{2-14}
&1 &0.9183 &0.9323 &0.9478 & 0& 0.8&11 &10 &46.8 &208 &10 &46 &197 \\
\cline{2-14}
&2 &0.8483 &0.8711 &0.9034 & 0&0.2 &3 &4 &18.5 &55 & 4&18.3 &55 \\
\hline
\hline
\multirow{5}{*}{30}&0.05  &0.9966 &0.9971 &0.9975 &0 &1.5 &13 &537 &1846 & 6194&537 &1844.4 &6186 \\
\cline{2-14}
&0.1 & 0.9933& 0.9941&0.9952 &0 &4.2 &41 &167 &1584.8 &24127 &166 &1580.6 &24094 \\
\cline{2-14}
&0.5 & 0.9655&0.9711 &0.9763 &0 &1.3 &11 &28 &127.5 &606 &28 &126.2 &596 \\
\cline{2-14}
&1 & 0.9345&0.9438 &0.9543 &0 &0.9 &9 &18 &59.4 & 359& 17& 58.5&350 \\
\cline{2-14}
&2 & 0.8758& 0.8928&0.9094 &0 &0.2 &2 &5 &22.9 &85 &5 &22.7 &85 \\
\hline
\hline
\hline
\end{tabular}
}
\caption{Benchmark results for a convex quadratic objective function}
\label{convextable}
\end{table}

\subsubsection{Comments About the Convex Case} 
While our running examples introduced at Subsection~\ref{runningex},  $\kopt$ could cover the whole set of natural integers (using for example $\gamma=\exp(-1/n)$). With those randomly generated matrices, even if the spectral radius is relatively close to 1, the integer $\kopt$ attains 41 (at $d=30$ and for $\varepsilon=0.1$) in the convex case. Most of generated problems have an integer $\kopt$ equal to 0. However, the maximal and the average integer $\kopt$ slightly decrease with respect to the spectral radius. This fact is independent of the dimension of the system. The maximal and the average $\kopt$ seem to be low-sensitive with respect to the dimension of the system and slightly increase when the dimension grows. 

The integer $\klyap_\nu$ is very sensitive to the dimension and the spectral radius of the system matrix. This is highlighted by the column on the minimal values of $\klyap_\nu$. This is not so surprising since $\klyap_\nu$ depends increasingly on $\mu(P)$ and $\norm{A}_P$. The value $\mu(P)$ is lower-bounded by $\lambda_{\rm min}(P)\max_{x\in\xin} \norm{x}_2^2$ and the construction of $\xin$ at Eq.~\eqref{initeq} implies that $ \max_{x\in\xin} \norm{x}_2^2$ grows with the dimension. On the other hand, $\norm{A}_P$ is lower-bounded by the spectral radius of $A$ and thus tends to grow when the spectral radius grows.  

The values at the last column about the differences between $\klyap_\nu$ and $\kopt$ are quite close to the ones of column $\klyap_\nu$. This is explained by the fact that $\kopt$ is mostly null.
We remark that for $d=3$ (and $\varepsilon=0.5,1,2$) and $d=5$ ($\varepsilon=2$), $\klyap_\nu$ is optimal since the difference between $\klyap_\nu$ and $\kopt$ is equal to one. 
\begin{table}[h!]
\centering
{\footnotesize
\begin{tabular}{|>{\centering }p{0.25cm}|>{\centering }p{0.5cm}|>{\centering }p{0.7cm}>{\centering }p{0.7cm}>{\centering }p{0.7cm}|>{\centering }p{0.2cm}>{\centering }p{0.3cm}>{\centering }p{0.3cm}|>{\centering }p{0.4cm}>{\centering }p{0.55cm}>{\centering }p{0.75cm}|>{\centering }p{0.5cm}>{\centering }p{0.6cm}>{\centering\arraybackslash }p{0.75cm}|}
\hline
\hline
\hline
\multicolumn{14}{|c|}{Concave}\\
\hline
d & $\varepsilon$ & \multicolumn{3}{c|}{Spec A}& \multicolumn{3}{c|}{$\kopt$} & \multicolumn{3}{c|}{$\klyap_\nu$}  &\multicolumn{3}{c|}{$\klyap_\nu-\kopt$}\\
\hline
& &Mi & Av & Mx & Mi & Av & Mx & Mi & Av & Mx &  Mi & Av & Mx \\
\hline
\hline
\multirow{5}{*}{3} & 0.05 &0.9574 & 0.9888& 0.9946
&0 &15.7 &106 & 127& 810.4&7610 & 85&794.7 &7594 \\
\cline{2-14}
&0.1 &0.9532& 0.9789& 0.9892& 0& 7.1& 39& 97 & 370.1&4624 &96 & 363&4590  \\
\cline{2-14}
&0.5 &0.8252 & 0.904&0.944 & 0&2.5 &15
&15 &78.8 & 638& 15 &76.3 & 623\\
\cline{2-14}
&1 &0.6649& 0.8313& 0.9011& 0& 1.2& 6 & 10&27.3 & 115& 10& 26.1&111  \\
\cline{2-14}
&2 &0.5338 & 0.7056&0.8089 & 0&0.4&3
&5 &11.9 & 34& 5 &11.4 & 34\\
\hline
\hline
\multirow{5}{*}{5}&0.05 &0.9814 &0.9922 &0.995 &0 &16.1 &176 &287 &1371.5 &10804 &281 & 1355.4&10746\\
\cline{2-14}
&0.1 &0.9742&0.9851 &0.9917 &0&9.1 &49 &155 &729 &17238 & 140&720.5 &17219 \\
\cline{2-14}
&0.5 & 0.871& 0.9285& 0.9577&0 &2.5 &17 &30 &110.3 &901 &28 &107.7 &888\\
\cline{2-14}
&1 &0.8032 &0.8665 &0.9238 &0 &1 &8 &17 &49.4 &162 &15 &48.4 &162 \\
\cline{2-14}
&2 & 0.6684&0.7627 & 0.8358& 0&0.6 &6 &10 &20.9 &110 &10 &20.2 &108\\
\hline
\hline
\multirow{5}{*}{10} & 0.05 &0.9927& 0.9949& 0.9962& 0& 19.4&234 &730 &5276.1 &250973 &730 &5256.6 &250839 \\
\cline{2-14}
&0.1 & 0.9861& 0.9897& 0.9924&0 &11.5 &77 &288 &1348.6 & 21116&286 &1337.1 &21039 \\
\cline{2-14}
&0.5 & 0.9363& 0.9492&0.9638 &0 &2.8 &15 &76 &223 &1835 &76 &220.5 &1829 \\
\cline{2-14}
&1 & 0.8677&0.9046 &0.9335 &0 &1.9 &9 &32 &106.1 & 684&32 &104.2 &682 \\
\cline{2-14}
&2 & 0.7707& 0.8274& 0.8723&0 &1.3 &6 &18 &37.3 &156 &16 &36 &155 \\
\hline
\hline
\multirow{5}{*}{20}& 0.05 &0.9956 & 0.9964&0.9971 &0 &21.7 &185 &1296 &4923.7 &124032 &1280 &4902 &123926 \\
\cline{2-14}
&0.1 &0.9905 &0.9927 &0.9945 &0 &12.2 &130 &603 &2367.8 &21960 &603 &2355.6 &21960 \\
\cline{2-14}
&0.5 & 0.9552& 0.9643& 0.9714&0 &4.7 &28 &115 &551 &5002 &115 &546.4 &5001 \\
\cline{2-14}
&1 &0.9167 &0.9315 &0.9461 &0 &2.6 &11 &67 &160.4 &756 &65 &157.8 &753 \\
\cline{2-14}
&2 & 0.8486& 0.8704&0.8942 &0 &1.6 &7 &35 &72 &373 &35 &70.4 &366 \\
\hline
\hline
\multirow{5}{*}{30} &0.05 &0.9967 &0.997 &0.9976 &0 &31 &335 &1643 &6007.8 &101884 &1635 &5976.8 &101821 \\
\cline{2-14}
& 0.1 & 0.993&0.994 &0.9949 &0 &12.2&73 &822 &4732.9 &128016 &820 &4720 &127973 \\
\cline{2-14}
&0.5 & 0.967& 0.9712& 0.9766& 0& 4.5& 43& 145&569.7 &6466 &143 &565.1 &6465 \\
\cline{2-14}
&1 & 0.9365& 0.9435& 0.9533& 0&3.7 &19 &76 &210.6 &856 &76 &206.9 &850 \\
\cline{2-14}
&2 & 0.8815& 0.893&0.9078 &0 &1.8 &9 &47 &90 &273 &47 &88.1 &273 \\
\hline
\hline
\hline
\end{tabular}
}
\caption{Benchmark results for a concave quadratic objective function}
\label{concavetable}
\end{table}

\subsubsection{Comments About the Concave Case} 
The integers $\kopt$ are quite higher compared to the convex case. We achieve a maximal $\kopt$ equal to 335 ($d=30$ and $\varepsilon=0.05$).
This particularity could be explained by the specific relation between $Q$ and $q$. The sensibility with respect to the dimension is more evident. The average values of $\kopt$ increase when the dimension increases  even if the maximal $\kopt$ for $d=10$ and $\varepsilon=0.05$ is higher than the maximal $\kopt$ for $d=30$ and $\varepsilon=0.05$. Finally, as in the convex case, most of $\kopt$ are equal to 0. 

We find the same sensibility of $\klyap_\nu$ with respect to the dimension and the spectral radius of the system matrix. Even if the formula changes since $Q\preceq 0$, $\klyap_\nu$ still depends increasingly on $\mu(P)$ and $\norm{A}_P$.

The values at the last column about the differences between $\klyap_\nu$ and $\kopt$ keep close to the ones of the column $\klyap_\nu$. However, the optimality of our formula is not attained in the context of those concave objective functions (the integers of the column Mi are strictly greater than 1).

\begin{table}[h!]
\centering
{\footnotesize
\begin{tabular}{|>{\centering }p{0.25cm}|>{\centering }p{0.5cm}|>{\centering }p{0.7cm}>{\centering }p{0.7cm}>{\centering }p{0.7cm}|>{\centering }p{0.2cm}>{\centering }p{0.3cm}>{\centering }p{0.3cm}|>{\centering }p{0.4cm}>{\centering }p{0.55cm}>{\centering }p{0.75cm}|>{\centering }p{0.5cm}>{\centering }p{0.6cm}>{\centering\arraybackslash }p{0.75cm}|}
\hline
\hline
\hline
\multicolumn{14}{|c|}{Linear}\\
\hline
d & $\varepsilon$ & \multicolumn{3}{|c|}{Spec A}& \multicolumn{3}{|c|}{$\kopt$} & \multicolumn{3}{|c|}{$\klyap_\nu$}  &\multicolumn{3}{|c|}{$\klyap_\nu-\kopt$}\\
\hline
& &Mn & Av & Mx & Mn & Av & Mx & Mn & Av & Mx &  Mn & Av & Mx \\
\hline
\hline
\multirow{5}{*}{3} & 0.05 &0.9729 & 0.9885 & 0.9944 & 0 & 2.4 & 25 & 6 & 384.1 & 8363 & 5 & 381.7 & 8362 \\
\cline{2-14}
& 0.1 &0.9564 & 0.9783 & 0.9899& 0 & 1.9 & 14& 4 & 119 & 949 &3 & 117.1 & 944  \\
\cline{2-14}
& 0.5 &0.7796 & 0.8996 & 0.9493 & 0 & 1.1 & 6 &1 & 25.9 & 645& 1 & 24.8 & 645\\
\cline{2-14}
& 1 &0.6697 & 0.827 & 0.8977 & 0 & 0.7 & 4 &1 & 8.3 & 25& 1 & 7.6 & 25\\
\cline{2-14}
& 2 &0.551 & 0.7047 & 0.8378 & 0 & 0.3 & 2 &1 & 3.9 & 22& 1 & 3.5 & 22\\
\hline
\hline
\multirow{5}{*}{5} & 0.05 &0.9887&0.9924 &0.9955 &0&3.4 &47 &30 &720.1 &12731 & 30&716.6 &12019 \\
\cline{2-14}
&0.1 &0.967 &0.9851 &0.9907 &0 &1.9 &17 &12&262.7 &6897 &12 & 260.7&6896\\
\cline{2-14}
&0.5 &0.8993&0.9294 &0.9595 &0&1 &7 &6 &40.4 &435 & 4&39.3 &434 \\
\cline{2-14}
&1 & 0.8151& 0.8704& 0.9161&0 &0.7 &5 &2 &15.4 &64 &2 &14.6 &59\\
\cline{2-14}
&2 &0.6847 &0.7619 &0.8313 &0 &0.5 &6 &1 &8.2 &32 &1 &7.6 &30 \\
\hline
\hline
\multirow{5}{*}{10}&0.05 & 0.9929& 0.9948& 0.9965& 0& 3.5&28 &364 &3797.8 &210030 &364 &3794.2 &210002 \\
\cline{2-14}
&0.1 & 0.9874& 0.9899& 0.9931&0 &2.8 &14 &91 &686.6 & 5941&91 &683.8 &5941 \\
\cline{2-14}
&0.5 & 0.9337& 0.9497&0.9632 &0 &1.9 &14 &30 &112.2 &1348 &29 &110.3 &1334 \\
\cline{2-14}
&1 & 0.874&0.9064 &0.9336 &0 &1.2 &7 &9 &39.5 & 177&9 &38.3 &176 \\
\cline{2-14}
&2 & 0.7755& 0.8278& 0.8765&0 &0.8 &6 &3 &17.7 &65 &3 &16.9 &62\\
\hline
\hline
\multirow{5}{*}{20}&0.05 &0.9955 & 0.9964&0.9971 &0 &3 &18 &628 &2183.4 &23913 &627 &2180.4 &23895 \\
\cline{2-14}
& 0.1 &0.991 &0.9927 &0.9945 &0 &2.9 &19 &269 &914.6 &4206 &269 &911.6 &4203 \\
\cline{2-14}
&0.5 & 0.9569& 0.9647& 0.971&0 &2.1 &11 &53 &166.1 &958 &53 &164 &950 \\
\cline{2-14}
& 1 &0.9169 &0.9309 &0.9448&0 &2 &8 &28 &80.1 &329 &25 &78 &325 \\
\cline{2-14}
& 2 & 0.8482& 0.8723&0.9023 &0 &1 &6 &8 &33.4 &115 &8 &32.3 &115 \\
\hline
\hline
 \multirow{5}{*}{30}& 0.05 &0.9965 &0.997 &0.9975 &0 &4 &37 &873 &4949.9 &85580 &872 &4945.9 &85575 \\
\hline
&0.1 & 0.9932&0.9941 &0.9952 &0 &3.2&23 &378 &2028 &38066 &378 &2024.8 &38054 \\
\hline
&0.5& 0.9667& 0.9711& 0.9765& 0& 2.1& 12& 62&270 &2333 &62 &267.8 &2323 \\
\hline
&1 & 0.9365& 0.9437& 0.9542& 0&2 &17 &84 &112.9 &351 &34 &110.9 &346 \\
\hline
&2 & 0.8812& 0.8935&0.9105 &0 &1.1 &5 &15 &48.2 &122 &15 &47 &121 \\
\hline
\hline
\hline
\end{tabular}
}
\caption{Benchmark results for a linear objective function}
\label{lineartable}
\end{table}

\subsubsection{Comments About the Linear Case} 
The results for the linear case are quite similar to the convex case. The main differences concern the dependency of $\kopt$ with respect to the dimension and the integers $\klyap_\nu$ for which the column "Mx" contains bigger integers. The average values of $\kopt$ seem to be less sensitive compared to the convex case whereas maximal values of $\kopt$ do not follow a monotonic progression.

\subsubsection{Graphical Summary of Benchmark Tests}
To obtain an user friendly view, we propose two graphics to sum up the results of these benchmarks. At left, we plot, in the logarithmic scale, the averages of the averages of $\kopt$ and $\klyap_\nu$ with respect to the dimension of the system for each type of problem (convex, concave, linear). The graphic at right represents, in the logarithmic scale, the averages of the averages of $\kopt$ and $\klyap_\nu$ with respect to the perturbation $\varepsilon$ applied to the spectral radius of the randomly generated matrix for each type of problems (convex, concave, linear).

For both graphics, the dotted line represents the averages for $\klyap_\nu$ whereas the plain lines represents the averages for $\kopt$. 
\begin{center}
\begin{figure}[h!]
    \begin{tikzpicture}[scale=0.9]
    \begin{scope}[xshift=0\textwidth]
    \begin{axis}[
    axis lines=left,
    axis x line =bottom,
    ymax=10,
    ymin=-2.5,
    xmax=33,
    xmin=0,
    minor tick num=5,
    x tick label style={
    below
    },
    x label style={at={(1.02,-0.01)},anchor=south},
    y label style={at={(0,1.15)},rotate=-90,anchor=north},
    xlabel = $d$,
    ylabel = {$k_{\rm opt}, K_\nu^{\rm lyap}$},
    ]
    \addplot[thick,red] coordinates {(3,{ln(1.06)})  (5,{ln(1.14)})  (10,{ln(1.66)})  (20,{ln(1.84)})  (30,{ln(1.62)})};
     \label{cvkp}
     \addplot[thick,dashed,red] coordinates{(3,{ln(67.44)}) (5,{ln(140.8)}) (10,{ln(355.88)}) (20,{ln(725.24)})  (30,{ln(728.12)})}; 
      \label{cvkl}
 \addplot[thick,blue]  coordinates{(3,{ln(5.38)}) (5,{ln(5.86)}) (10,{ln(7.38)}) (20,{ln(8.56)}) (30,{ln(10.64)})};
 \label{cckp}
 \addplot[thick,dashed,blue]  coordinates{(3,{ln(259.7)}) (5,{ln(456.22)}) (10,{ln(1398.2)}) (20,{ln(1615)}) (30,{ln(2322.2)})};
\label{cckl}
 \addplot[thick,black]  coordinates{(3,{ln(1.28)})  (5,{ln(1.5)}) (10,{ln(2.04)}) (20,{ln(2.2)}) (30,{ln(2.48)})};
 \label{likp}
    \addplot[thick,dashed,black]  coordinates{(3,{ln(108.24)}) (5,{ln(209.36)}) (10,{ln(930.76)}) (20,{ln(675.52)}) (30,{ln(1481.8)})};
    \label{likl}
\end{axis}
\end{scope}
\begin{scope}[xshift=0.6\textwidth]
 \begin{axis}[
    axis lines=left,
    axis x line =bottom,
    ymax=10,
    ymin=-2.5,
    xmax=2.2,
    xmin=0,
    minor tick num=5,
    x tick label style={
    below
    },
    x label style={at={(1.02,-0.01)},anchor=south},
    y label style={at={(0,1.15)},rotate=-90,anchor=north},
    xlabel = $\varepsilon$,
    ylabel = {$k_{\rm opt}, K_\nu^{\rm lyap}$}
    ]
    \addplot[thick,black]  coordinates{(0.05,{ln(3.26)}) (0.1,{ln(2.54)}) (0.5,{ln(1.5)}) (1,{ln(1.32)}) (2,{ln(0.74)})};
    \addplot[thick,dashed,black]  coordinates{(0.05,{ln(2407.1)}) (0.1,{ln(802.18)}) (0.5,{ln(122.92)}) (1,{ln(51.24)}) (2,{ln(22.28)})};
    \addplot[thick,blue]  coordinates{(0.05,{ln(20.78)}) (0.1,{ln(10.42)}) (0.5,{ln(3.4)}) (1,{ln(2.08)}) (2,{ln(1.14)})};
    \addplot[thick,dashed,blue]  coordinates{(0.05,{ln(3677.4)}) (0.1,{ln(1909.7)}) (0.5,{ln(306.56)}) (1,{ln(110.76)}) (2,{ln(46.42)})};
    \addplot[thick,red]  coordinates{(0.05,{ln(2.42)}) (0.1,{ln(2.74)}) (0.5,{ln(1.14)}) (1,{ln(0.8)}) (2,{ln(0.22)})};
    \addplot[thick,dashed,red]  coordinates{(0.05,{ln(1142.5)}) (0.1,{ln(753.14)}) (0.5,{ln(74.74)}) (1,{ln(32.02)}) (2,{ln(12.14)})};
\end{axis}
\end{scope}
 \matrix [
        draw,
        matrix of nodes,
        anchor=north,
        yshift=-2.5mm,
        node font=\small,
    ] at (8,-1) {
        Concave        & Convex     & Linear \\
        \ref{cckp}  &  \ref{cvkp}  &  \ref{likp} \\
        \ref{cckl} &   \ref{cvkl}&  \ref{likl}  \\
    };
\end{tikzpicture}
\caption{Graphical summary of benchmark results}
\end{figure}
\end{center}
On the left, for $\kopt$, we observe a slight grow for the concave and convex cases whereas the monotony is not significant for the linear case. For $\klyap_\nu$, the same observations holds. However, the linear case is still between the convex and the concave cases. 
On the right, both $\kopt$ and $\klyap_\nu$ decrease as the deviation increases. In other words, when the spectral radius decreases the integers decrease. We find the same order on the type of problem as at left : the linear case is still between convex and concave type.

\section{Conclusion and Future Works}
\label{conclusion}
We propose a computable method to solve a particular maximization problem where the objective function is quadratic (or linear) and where the constraints set is the reachable set of a stable affine dynamical system. This problem can be viewed as to search the maximal term of a upper-bounded real sequence. Our solution is to compare the first $N$ terms of the sequence. The computation of the integer $N$ guarantees the fact that the maximal term has been found within the $N$ first terms of the sequence. The computation of $N$ derives from a formula based on some classical linear algebra tools such as maximal eigenvalues of symmetric matrices, matrix norms and discrete Lyapunov equation solutions. The formula proposed is flexible and can be optimized in order to minimize the number of comparisons. 

The work proposed is an extension of the case where the system matrix was only diagonalizable. Here the system matrix is just stable. We illustrate our approach on a family of nondiagonalizable matrices of dimension two. In this special case, we are able to reduce the number of comparaisons from a dichotomy method. In the general case, the problem is a nonlinear semidefinite program with a nondifferentiable objective function. Future works would concern the development of a resolution of this particular nonlinear semidefinite program. 
\bibliographystyle{tfnlm}
\bibliography{versionfinalbib}
\newpage

\section*{Appendix}

\subsection*{Inequalities of Example~\ref{complicated}}
The final inequality to prove is the following : for all $k\in\nn$ and all $B\in (1,\gamma^{-2})$:
\[
(1+k)^2\leq \left(1+\dfrac{B}{(B-1)^2}\right)B^k
\]
As $\gamma$ browse $(0,1)$, this is the same as proving the inequality for $B\in(1,+\infty)$.

Applying the natural logarithm, we have to prove for all $k\in\nn$, for all $B\in(1,+\infty)$:
\[
2\ln(1+k)\leq k\ln B+\ln\left(1+\dfrac{B}{(B-1)^2}\right)
\]
First, the sequence, $(\ln((1+k)^2)-k\ln(B))_k$ is decreasing if and only if for all $k\in\nn$, $2\ln(1+(1+k)^{-1})\leq \ln B$ which is equivalent to $\ln(4)\leq \ln B$ which is the same as $4\leq B$. We conclude by the fact that for $k=0$, the inequality holds. The inequality is proved for $B\geq 4$.

Now, we suppose that $B\in (1,4)$. The function $x\mapsto \ln(1+B(B-1)^{-2})+x\ln(B)-2\ln(1+x)$ is minimal at $x=(2/\ln(B))-1$ and thus it suffices to prove that the function
$B\mapsto \ln(1+B(B-1)^{-2})+2-\ln(B)-2\ln(2/\ln (B))$
is nonnegative. With simple calculus and the natural logarithm properties, this is the same as proving 
$e^2 ((B-1)^2+B)(\ln(B))^2- 4B(B-1)^2\geq 0$.
We have 
$e^2 ((B-1)^2+B)(\ln(B))^2-4B(B-1)^2
= (e\ln(B)\sqrt{(B-1)^2+B}-2\sqrt{B}(B-1))(e\ln(B)\sqrt{(B-1)^2+B}+2\sqrt{B}(B-1))$.
Since $B>1$,  $e\ln(B)\sqrt{(B-1)^2+B}+2\sqrt{B}(B-1)>0$ and we must prove that $e\ln(B)\sqrt{(B-1)^2+B}-2\sqrt{B}(B-1)\geq 0$. We have: $e\ln(B)\sqrt{(B-1)^2+B}-2\sqrt{B}(B-1)
\geq e\ln(B)\sqrt{(B-1)^2+B-1}-2\sqrt{B}(B-1)\geq \sqrt{B-1}\sqrt{B}(e\ln(B)-2\sqrt{B-1})$.
It suffices to prove that $e\ln(B)-2\sqrt{B-1}\geq 0$ .
We write $X=\sqrt{B-1}\in (0,\sqrt{3})$ and thus $B=1+X^2$. We study
$X\mapsto \varphi(X)=e\ln(1+X^2)-2X$ on $(0,\sqrt{3})$. The derivative of $\varphi$ is $(2eX-2(1+X^2)/(1+X^2)$. The polynomial
$-X^2+eX-1$ has two real roots $(e-\sqrt{e^2-4})/2$ and $(e+\sqrt{e^2-4})/2\geq \sqrt{3}$. Hence, we have the following variations for $\varphi$:
\begin{center}
\begin{tikzpicture}[scale=0.8]
   \tkzTabInit[deltacl = 1.2]{$X$ / 1 , $\varphi'(X)$ / 1, $\varphi(X)$ / 2}{$0$, $(e-\sqrt{e^2-4})/2$, $\alpha$ , $\sqrt{3}$}
   \tkzTabLine{,-,z ,+ ,, +}
   \tkzTabVar{+/0,-/ , R/, +/ ,R/, $e\ln(4)-2\sqrt{3}$}
   \tkzTabIma{2}{4}{3}{$0$}
\end{tikzpicture}
\end{center}
With a dichotomy, we can prove that $\alpha\in (1.1637,1.1638)$. Thus, we conclude that $\varphi$ is positive between $1.2$ and $\sqrt{3}$. This means that $(e\ln(B)\sqrt{(B-1)^2+B}-2\sqrt{B}(B-1))$ is nonnegative for $B\in [2.44,4)$.

To complete the proof, we have to show that $e\ln(B)\sqrt{(B-1)^2+B}-2\sqrt{B}(B-1)$ is nonnegative for $B\in (1,2.44)$. We use a sharp lower bound for the logarithm~\cite{b0750f61-9854-3f22-8cec-fb4c23a27312} : for all $X>1$, $\ln(X)\geq 2(X-1)/(1+X)$. If we prove that for all $B\in(1,2.44]$, 
$2e(B-1)(B+1)^{-1}\sqrt{(B-1)^2+B}-2\sqrt{B}(B-1)\geq 0$
or equivalently for all $B\in(1,2.44]$, $2(B-1)(B+1)^{-1}(e\sqrt{(B-1)^2+B}-\sqrt{B}(B+1))\geq 0$, the result holds. This is the same as proving that:
$e\geq \sqrt{(B(B+1)^2)/((B-1)^2+B)}$.
Let us write for $B>1$, $\psi(B)=(B(B+1)^2)/(B^2-B+1)$. 
Then, for $B>1$: $\psi'(B)=(B+1)\left[B^3-3B^2+3B+1\right]/((B^2-B+1)^2)\geq (B+1)(B-1)^3/((B^2-B+1)^2)\geq 0$.
We conclude that $\psi$ is increasing on $(1,+\infty)$ and thus for all $B\in(1,2.44]$, $\psi(B)\leq \psi(2.44)<6.3972\leq e^2\simeq 7.3891$. This concludes the proof.

The second iinequality of Example~\ref{complicated} is :
\[
\nuopt<1+\gamma^2(1-\gamma^2)^{-2}
\]
If $0<\gamma\leq \exp(-1)$, $\nuopt=1$ and the strict inequality holds since $\gamma>0$. Now suppose that $\gamma >\exp(-1)$.  We have proved that the maximizer of  $\psi_\gamma:u\mapsto \gamma^{2u}(1+u)^2 $ is $u^*:=-1/\ln(\gamma)-1$. Hence, we have $\nuopt\leq (\exp(-2)(\gamma \ln(\gamma))^{-2}$. Thus $\nuopt<1+\gamma^2 (1-\gamma^2)^{-2}$ if $(\ln \gamma)^2(1-\gamma^2+\gamma^4)-\exp(-2)(\gamma^{-2}-2+\gamma^2)>0$. Now since $\gamma>\exp(-1)$, $(\ln\gamma)^2>1$ and it suffices to prove that $\gamma\mapsto g(\gamma):=(1-\gamma^2+\gamma^4)-\exp(-2)(\gamma^{-2}-2+\gamma^2)$ is strictly positive on $(\exp(-1),1)$. Since $g'(\gamma)=\gamma^{-3}(4\gamma^6-2(\exp(-2)+1)\gamma^4+2\exp(-2))$ and
the derivatives of $\gamma\mapsto \gamma^3 g'(\gamma)= \gamma^3(24 \gamma^2-8(\exp(-2)+1))$, we conclude that  $\gamma\mapsto \gamma^3 g'(\gamma)$ achieves its minimum at $\overline{\gamma}=\sqrt{3^{-1}(\exp(-2)+1)}$ for which $g'(\overline{\gamma})>0$. Finally, $g$ is strictly increasing and thus, for all $\gamma\in (\exp(-1),1)$, $g(\gamma)>g(\exp(-1))=\exp(-2)>0$. 

\subsection*{Proof of the Existence and Uniqueness of  A Minimizer and Convergence of the Dichotomy Method}

\begin{lemma}
\label{critical}
Let $g:\rr\mapsto \rr$ a $\mathcal{C}^1$ function. Suppose that $g$ has a critical point. If all critical points are strict local minimizers, then $g$ has a unique critical point which is the unique minimizer of $g$.
\end{lemma}
\begin{proof}
Let $\oB$ be a critical point and suppose that there exists another critical point. Wlog, we can suppose that $\oB$ is not the greatest and consider $C^*:=\min\{B'>\oB: h'(B')=0\}$. As $C^*$ is a critical point it is a strict local minimizer and there exist $B_0,C_0$ such that $\oB<B_0<C_0<C^*$, $h'(B_0)>0$ and $h'(C_0)<0$. By continuity of $h'$, there exists $C_1\in (B_0,C_0)$ such that $h'(C_1)=0$ which contradicts the minimality of $C^*$. This proves the uniqueness of the critical point.

Now, if there exists $C$ s.t. $h(C)<h(\oB)$. We can suppose that $C<\oB$. For some $\alpha>0$, $0<t<\alpha$ implies that $h(\oB-t)>h(\oB)>h(C)$ and $C<\oB-t$. Then there exists $D\in (C,\oB-t)$ such that $h(D)=h(\oB)$ and by Rolle's theorem, we conclude that $h'(C^*)=0$ for some $C^*\in (C,\oB-t)$. This is another critical point which is thus locally minimal. This contradicts the fact that $\oB$ is the only local minimizer.
\end{proof}

Recall that we define \[g_k:B\mapsto (\ln(B\gamma^2))^{-1}\ln\left(\dfrac{(B-1)^2\nu_k}{B^2-B+1}\right)\enspace.\] 
on $(1,\gamma^{-2})$.

The following assertions are true:
\begin{enumerate}
\item The nonlinear equation $g'_k(B)=0$ has a solution on $(1,\gamma^{-2})$;
\item For all $\overline{B}\in(1,\gamma^{-2})$ such that $g'_k(\overline{B})=0$, we have $g_k''(\overline{B})>0$. In other words, all critical points are local strict  minimizers. We conclude that there exists one critical point which is the unique minimizer of $g_k$; 
\item The dichotomy method converges to this unique minimizer.
\end{enumerate}
\begin{proof}
Let $B\in (1,\gamma^{-2})$. Let us write $f_k(B)=\frac{\nu_k(B-1)^2}{B^2-B+1}$. Hence, we have $g_k(B)=\ln(f_k(B))(\ln(\gamma^2 B))^{-1}$, $g'_k(B)=\left(f'_k(B)(f_k(B))^{-1}\ln(\gamma^2 B)-B^{-1}\ln(f_k(B))\right)(\ln(\gamma^2 B))^{-2}$ and $f'_k(B)=\nu_k(B^2-1)(B^2-B+1)^{-2}$. Now, we remark that $g'_k(B)=0$ if and only if $h_k(B):=f'_k(B)(f_k(B))^{-1}\ln(\gamma^2 B)-B^{-1}\ln(f_k(B))=0$. This implies that $h_k(B)$ tends to $-\gamma^2\ln(f_k(\gamma^{-2}))$ as $B$ tends to $\gamma^{-2}$ by continuity of $\ln \circ f_k$.
Writing $h_k(B)=(B-1)^{-1}\left(\frac{(B+1)\ln(\gamma^2 B)}{B^2-B+1}-\frac{(B-1)(\ln(\nu_k)+2\ln(B-1)-\ln(B^2-B+1))}{B}\right)$, we can show that $h_k(B)$ tends to $-\infty$ as $B$ tends to $1$. Note that $-\gamma^2 \ln(f_k(\gamma^{-2}))>0$. Indeed, we have $f_k(\gamma^{-2})<1$ which is a consequence of $\nu_k\leq \nuopt<1+\gamma^2(1-\gamma^2)^{-2}$ which is already proved at Example~\ref{complicated}. Those two limits and the continuity of $h_k$ imply the existence of some $\oB\in (1,\gamma^{-2})$ such that $h_k'(\oB)=0$.

We have proved that a critical point of $g_k$ exists. Now, we prove that any critical point of $g_k$ is a strict local minimizer i.e. if $g_k'(B)=0$ implies that $g_k''(B)>0$.
We have $g''_k(B)=\left((h_k(B))(\ln(\gamma^2 B))^{-1}\right)'=\left(h'_k(B)\ln(\gamma^2 B)^{2}-2 B^{-1} h_k(B) \ln(\gamma^2 B)\right)\left(\ln(\gamma^2 B)^{4}\right)^{-1}$. Then, when $g'_k(B)=0$ or equivalently $h_k(B)=0$, $g''_k(B)>0$ if and only if $h'_k(B)>0$. We also have:
\[
h'_k(B)=
\dfrac{(f''_k(B)\ln(\gamma^2 B)+B^{-1} f'_k(B))f_k(B)-(f'_k(B))^2\ln(\gamma^2 B)}{f_k(B)^2}-\dfrac{f'_k(B)}{Bf_k(B)}+\dfrac{\ln(f_k(B)}{B^2}
\]
and so:
\[
h'_k(B)=\dfrac{f''_k(B)\ln(\gamma^2 B)f_k(B)-(f'_k(B))^2\ln(\gamma^2 B)}{f_k(B)^2}+\dfrac{\ln(f_k(B)}{B^2}\enspace .
\]
Note that $h_k(B)=0$ is equivalent to $\dfrac{\ln(f_k(B))}{B^2}=\dfrac{f'_k(B)}{Bf_k(B)}\ln(\gamma^2 B)$. Finally, for all $\oB$ such that $h_k(\oB)=0$ we have:
\[
h'_k(\oB)=\dfrac{f''_k(B)\ln(\gamma^2 B)f_k(B)-(f'_k(B))^2\ln(\gamma^2 B)}{f_k(B)^2}+\dfrac{f'_k(B)}{Bf_k(B)}\ln(\gamma^2 B)
\]
and then:
\[
h'_k(\oB)=\dfrac{\ln(\gamma^2 \oB)}{\oB f_k(\oB)^2}\biggl(\oB f''_k(\oB)f_k(\oB)-\oB (f'_k(\oB))^2+f_k(\oB)f'_k(\oB)\biggr)\enspace .
\]
Recall that $f'_k(B)=\nu_k\dfrac{B^2-1}{(B^2-B+1)^2}$ and thus 
$f_k''(B)=-2\nu_k\dfrac{B^3-3B+1}{(B^2-B+1)^3}$.
We obtain:
\[
\oB f''_k(\oB)f_k(\oB)-\oB (f'_k(\oB))^2+f_k(\oB)f'_k(\oB)
=\dfrac{\nu_k^2 (\oB-1)^2}{(\oB^2-\oB+1)^4}\left(-\oB^4-2\oB^3+4\oB^2-2\oB-1\right)
\]
Since $\oB>1$, we have finally:
\[
-\oB^4-2\oB^3+4\oB^2-2\oB-1=-(\oB^2-2)^2+3-2\oB^3-2\oB\leq 3-2\oB^3-2\oB<-1
\]
We conclude that $\oB f''_k(\oB)f_k(\oB)-\oB (f'_k(\oB))^2+f_k(\oB)f'_k(\oB)$ is negative and since $\oB\in (1,\gamma^{-2})$, $\ln(\gamma^2\oB)<0$, this leads to $h'_k(\oB)>0$ or equivalently $g_k''(\oB)>0$ which proves that $\oB$ is a strict local minimizer of $g_k$. This proves that all critical points are strict local minimizers. We conclude from Lemma~\ref{critical}. 
\end{proof}

\subsection*{Initial Polytope Generation for Benchmarks}
Recall that we use $\xin=\conv(\{x_k,k=1,\ldots,d+1\})$
where the family $\{x_k\}_{k=1,\ldots,d+1}$ is defined by :
\[
\begin{array}{lr}
x_{1,i}= -1, \forall\, i=1,\ldots d,& x_{2,i}=\left\{\begin{array}{cr} 1 & \text{ if } i=1\\ 0 & \text{ if } 2\leq i\leq d\end{array}\right.,\\
\\
 x_{k,i}=\left\{\begin{array}{cr} -1 & \text{ if } 1\leq i\leq k-2\\ 1 & \text{if } i=k-1 \\ 0 & \text{ if } k\leq i\leq d\end{array}\right. \text{ and } &  x_{d+1,i}=\left\{\begin{array}{cr} -1 & \text{ if } 1\leq i\leq d-1\\ 1 & \text{if } i=d \end{array}\right.
 \end{array}
\]
\begin{proposition}
The following statements hold:
\begin{enumerate}
\item The family $\{x_k\}_{k=1,\ldots,d+1}$ is affinely independant;
\item Let us define the matrix $F$ such that for all $i,j\in \{1,\ldots d\}$:
\[
F_{i,j}=\left\{
\begin{array}{lr} 
2^{j-i-1} & \text{ if } i>j\\
-\frac{1}{2} & \text{ if } i=j\\
0 & \text{ if }j > i
\end{array}
\right.
\]
Then :
\[
y\in\conv\left(\left\{x_k,k=1,\ldots,d+1\right\}\right)\iff \left\{\begin{array}{l}\displaystyle{F_i y\leq 2^{-i}},\ \forall\, i=1,\ldots,d\\ 
\\
\displaystyle{\sum_{j=1}^d} 2^{j-1} y_j\leq 1\end{array}\right.
\]
\item $0\in \inte\conv\left(\left\{x_k,k=1,\ldots,d+1\right\}\right)$.
\end{enumerate}
\end{proposition}
\begin{proof}
Let us define the matrix $C$ such that for all $i,j\in\{1,\ldots,d\}$:
\[
C_{ij}=\left\{
\begin{array}{lr}
-1 & \text{ if } i < j\\
1 & \text{ if } i=j\\
0 & \text{ if } i > j
\end{array}
\right.
\]
The matrix $C$ is invertible and for all $i,j\in\{1,\ldots,d\}$:
\[
C_{ij}^{-1}=\left\{
\begin{array}{lr}
2^{j-i-1} & \text{ if } i < j\\
1 & \text{ if } i=j\\
0 & \text{ if } i > j
\end{array}
\right.
\]
We also use the notation $\mm$ for the vectors for which the coordinates are equal to 1.

{\itshape 1}. The matrix for which the column are the vectors of $\{x_k-x_1\}_{k=2,\ldots,d+1}$ is the matrix $D$ s.t. $D_{ij}=C_{ij}+1$ for all $i,j\in \{1,\ldots,d\}$. The matrix  $D$ is lower triangular with diagonal elements equal to 2. Hence, its determinant is equal to $2^d$ and the family $\{x_k-x_1\}_{k=2,\ldots,d+1}$ is a basis of $\rd$ and the conclusion follows.

{\itshape 2}. Let $y\in\conv\left(\left\{x_k,k=1,\ldots,d+1\right\}\right)$. Then, there exist a nonnegative scalar $\mu$ and a nonnegative vector $\lambda=(\lambda_1,\ldots,\lambda_d)$ such that $y=-\mu \mm + C\lambda$ and $\mu+\mm^\intercal \lambda=1$. This is equivalent to $C^{-1}(y+\mu \mm)=\lambda$ which is the same as  $C^{-1}(y+\mu \mm)=\lambda$ and $\mu =(1+\mm^\intercal C^{-1}\mm)^{-1}(1-\mm^\intercal C^{-1}y)$. Now we can remove the coefficient $\mu$ and $\lambda$ to get:
\[
y\in\conv\left(\left\{x_k,k=1,\ldots,d+1\right\}\right)\iff
\left\{
\begin{array}{l}
C^{-1}\left(y+\dfrac{1-\mm^\intercal C^{-1}y}{1+\mm^\intercal C^{-1}\mm}\mm\right)\geq 0\\
\\
\dfrac{1-\mm^\intercal C^{-1}y}{1+\mm^\intercal C^{-1}\mm}\geq 0
\end{array}
\right.
\]
Note that $C$ and $C^{-1}$ have only one eigenvalue equal to 1 and thus are positive definite and $1+\mm^\intercal C^{-1}\mm$ is strictly positive. Hence $(1+\mm^\intercal C^{-1}\mm)^{-1}(1-\mm^\intercal C^{-1}y)$ is nonnegative if and only if $(1-\mm^\intercal C^{-1}y)$ is nonnegative. By simple calculus, we get $1+\mm^\intercal C^{-1}\mm=2^d$. Moreover, as $\mm^\intercal C^{-1}y$ is a scalar $\mm^\intercal C^{-1}y\mm= \mm\mm^\intercal C^{-1}y$ and thus:
\[
\begin{array}{ll}
\displaystyle{C^{-1}\left(y+\dfrac{1-\mm^\intercal C^{-1}y}{1+\mm^\intercal C^{-1}\mm}\mm\right)}&=\displaystyle{C^{-1}y+2^{-d}(C^{-1}\mm-C^{-1}\mm\mm^\intercal C^{-1}y)}\\
&=
\displaystyle{(C^{-1}-2^{-d}C^{-1}\mm\mm^\intercal C^{-1})y+2^{-d}C^{-1}\mm}
\end{array}
\]
Again, some simple calculus lead to, for all $i=1,\ldots,d$, $(C^{-1}\mm)_i=2^{d-i}$ and $(\mm^\intercal C^{-1})_i=2^{i-1}$ and for all $i,j\in\{1,\ldots,d\}$, $(C^{-1}\mm\mm^\intercal C^{-1})_{ij}=2^{d+j-i-1}$. We conclude by setting $F=2^{-d}C^{-1}\mm\mm^\intercal C^{-1}-C^{-1}$.

{\itshape 3}. We have $F_i 0 =0 < 2^{-i}$ for all $i=1,\ldots,d$ and $\sum_{j=1}^d=2^{j-1} 0=0<1$. We conclude that $0\in \inte\conv\left(\left\{x_k,k=1,\ldots,d+1\right\}\right)$. 
\end{proof}
\end{document}